# An ultraweak DPG method for viscoelastic fluids

B. Keith[a], P. Knechtges[b], N.V. Roberts[c,1], S. Elgeti[b], M. Behr[b], L. Demkowicz[a]

[a]*The Institute for Computational Engineering and Sciences, The University of Texas at Austin, Austin, Texas*
[b]*Chair for Computational Analysis of Technical Systems, RWTH Aachen University, Aachen, Germany*
[c]*Center for Computing Research, Sandia National Laboratories, Albuquerque, New Mexico*

## Abstract

We explore a vexing benchmark problem for viscoelastic fluid flows with the *discontinuous Petrov-Galerkin* (DPG) finite element method of Demkowicz and Gopalakrishnan [1, 2]. In our analysis, we develop an intrinsic *a posteriori* error indicator which we use for adaptive mesh generation. The DPG method is useful for the problem we consider because the method is inherently stable—requiring no stabilization of the linearized discretization in order to handle the advective terms in the model. Because stabilization is a pressing issue in these models, this happens to become a very useful property of the method which simplifies our analysis. This built-in stability at all length scales and the *a posteriori* error indicator additionally allows for the generation of parameter-specific meshes starting from a common coarse initial mesh. A DPG discretization always produces a symmetric positive definite stiffness matrix. This feature allows us to use the most efficient direct solvers for all of our computations. We use the Camellia finite element software package [3, 4] for all of our analysis.

## 1. Introduction

Viscoelastic fluids and the models which are used to predict their behavior are common in engineering mechanics and industry. Typical examples of such fluids are blood and polymer melts and interest for modeling them is often found in biomedical engineering [5] and plastics manufacturing [6]. Given their significance, it is important to note that the modeling of these fluids is challenging from both the computational as well as the physical perspective.

From the computational perspective, challenges manifest in three important ways:

1. Given the complexity of the equations, few analytical solutions are known, so that verification of numerical methods is relegated to benchmark problems. Ostensibly benign issues, discrepancies in numerical results and loss of convergence are common and contentious.
2. Stability is a prevalent concern due to the convective nature of the nonlinear viscoelastic constitutive laws. Given the discrepancies and eventual numerical failures reported in the literature, some suspect that stabilization may have introduced spurious effects in some discretizations.
3. The solution in the most controversial parameter ranges exhibit small-length-scale, high-contrast features. For instance, in the benchmark problem we considered, both boundary layers and internal layers develop in a stress variable. Local to these features, the values of the stress components rapidly change by orders of magnitude.

From the physical perspective—also due to the complexity of the equations—qualitative knowledge of solution features is limited. Indeed, it is not generally known in many benchmark problems upon which parameter ranges physical instabilities or transience is present. Therefore, there is contention in whether the discrepancies in benchmark results—and more importantly, the loss of convergence in numerical methods—are due to the ill-posedness of the very problems being discretized.

In an attempt to address the second and third computational issues mentioned above, we apply new finite element analysis to the *steady-state* confined cylinder benchmark problem with a method which guarantees numerical stability of each linearized problem *without stabilization terms* and involves a built-in *a posteriori* error estimator. These features will reduce the means for numerical error in the results and allow for adaptive resolution of the small-scale solution features via a sequence of refined meshes. Both of these aspects make our analysis notable in the context of the literature to date.

We will use the discontinuous Petrov-Galerkin (DPG) methodology of Demkowicz and Gopalakrishnan [7] for our investigation. Since it is rather new in the literature, we will contrast some features of DPG with those of more traditional stabilized finite element methods in the second part of this introduction. However, let us highlight here that the approach DPG takes is different from the log-conformation methods, as used in [8–13], in that it is stabilization for the linearized problem, whereas log-conformation is a non-linear reformulation of the problem. As such, the two approaches are independent and so it begs the question whether the combination of both methods may yield an even more improved numerical method. As similar point of

*Email address:* brendan@ices.utexas.edu (B. Keith)
[1]This research was supported in part by the Office of Science, U.S. Department of Energy, under Contract DE-AC02-06CH11357, through a postdoctoral appointment at Argonne National Laboratory. It was also supported in part by an appointment at Sandia National Laboratories, a multi-mission laboratory managed and operated by Sandia Corporation, a wholly owned subsidiary of Lockheed Martin Corporation, for the U.S. Department of Energys National Nuclear Security Administration under contract DE-AC04-94AL85000.

reasoning could also begin with the square-root method of [14].

We will now briefly expand upon the viscoelastic models that we will later analyze. The introduction closes with some comments on notation and an outline of the rest of the paper.

*1.1. Viscoelastic fluid models*

Introducing pressure, $p$, velocity, $\mathbf{u}$, and *solvent viscosity*, $\eta_\mathrm{S}$, a very common macroscopic description of (incompressible) viscoelastic fluids is given by the following constitutive law for the Cauchy stress [15]:

$$\boldsymbol{\sigma} = -p\,\mathbf{I} + \eta_\mathrm{S}\bigl(\boldsymbol{\nabla}\,\mathbf{u} + \boldsymbol{\nabla}\,\mathbf{u}^\mathsf{T}\bigr) + \mathbf{T}\,. \qquad (1.1)$$

Here, the non-Newtonian term, the *extra stress tensor*, $\mathbf{T}$, is governed by a relationship independent of kinematic conservation laws. For instance, the Giesekus model [16], with *mobility factor*, $\alpha \in [0, 1]$, describes the advection and decay of $\mathbf{T}$ along the streamlines of the fluid by

$$\mathbf{T} + \lambda\,\mathfrak{L}_\mathbf{u}\mathbf{T} + \alpha\frac{\lambda}{\eta_\mathrm{P}}\mathbf{T}^2 = \eta_\mathrm{P}\bigl(\boldsymbol{\nabla}\,\mathbf{u} + \boldsymbol{\nabla}\,\mathbf{u}^\mathsf{T}\bigr)\,. \qquad (1.2)$$

Here, $\lambda > 0$ is the *relaxation time*, $\eta_\mathrm{P}$ is the *polymeric viscosity*, and $\mathfrak{L}$ is the *Lie derivative* operator—in this situation, acting on $\mathbf{T}$ in the direction $\mathbf{u}$ and then often called the *upper-convected Maxwell derivative* [17]—viz.,

$$\mathfrak{L}_\mathbf{u}\mathbf{T} = \frac{\partial \mathbf{T}}{\partial t} + (\mathbf{u}\cdot\boldsymbol{\nabla})\mathbf{T} - (\boldsymbol{\nabla}\,\mathbf{u})\mathbf{T} - \mathbf{T}(\boldsymbol{\nabla}\,\mathbf{u})^\mathsf{T}\,. \qquad (1.3)$$

Notably, when $\alpha = 0$ the Giesekus model reduces to the Oldroyd-B model [18].

We will discretize both the Giesekus model and the Oldroyd-B model in this work. We note that there are many other common and similar viscoelastic models which we do not analyze. These include the White-Metzner model [19] and the Phan-Thien-Tanner (PTT) model [20] as well as various finitely extensible, nonlinearly elastic (FENE) models [21]. Even though the Oldroyd-B model suffers from a lack of finite extensibility, it is also the most well-studied of these models for the confined cylinder benchmark problem. Moreover, because most alternative models are closely related to Oldroyd-B, it is a natural place to start. So as to correlate our results with the literature even further, we additionally consider inertial effects in the Oldroyd-B model, and the (non-inertial) Giesekus model.

*1.2. DPG versus stabilized methods*

In the name "DPG", *Petrov*-Galerkin indicates that the trial and test spaces do not need to coincide. In fact, in this paper, they *will not coincide*. The DPG methodology is to assemble the test-space on-the-fly in such a way as to induce stability in each linearized problem. This is performed element-wise and made possible *because the test space is discontinuous*. Because of the local operations on the test space, the DPG methodology always produces a symmetric and positive definite stiffness matrix and is automatically stable with a stability constant predictably close to that of the infinite-dimensional problem.

*Bubnov*-Galerkin finite element methods use a discretization where both the trial and the test functions are drawn from identical function spaces. However, when there is a lack of symmetry in the equations—as is common in fluid flow problems—such choices of trial and test functions produce non-symmetric stiffness matrices. Generally speaking, non-symmetric matrices require more expensive linear solvers than symmetric matrices.

A further and perhaps more important hindrance is that Bubnov-Galerkin discretizations commonly produce unstable systems or systems with large stability constants. Rectifying this issue frequently involves introducing new terms and parameters into the discretization which are not present in the original equations. In some cases, the introduction of these new *stabilization* terms can be understood to be equivalent to modifications of the test space, so that many of these methods are also called Petrov-Galerkin methods.[2] With these methods, as the number of equations in the model grows, so does the size of the new parameter space and so also the difficulty of choosing suitable parameters. Methods which employ such strategies are frequently called *stabilized* methods and some of the most prominent and successful are the streamline upwind Petrov-Galerkin method [22] and the variational multiscale method [23].

In this paper, we apply the discontinuous Petrov-Galerkin finite element methodology [7] to an *ultraweak variational formulation* [1] of the steady version of the aforementioned viscoelastic fluid models. By an ultraweak variational formulation, we mean a variational formulation of a first-order PDE system wherein all of the derivatives have been moved from the trial functions and onto the test functions through integration by parts. Therefore, before we begin to discretize the problem, we must write the equations as a first-order system (see Section 3). Such formulations have been studied with DPG in both incompressible (Stokes and Navier-Stokes) [24–26] and compressible fluid flow problems (Navier-Stokes) [27, 28]. Demonstrated successes in these studies suggested that it would also perform well with viscoelastic fluid models. And, naturally, *without stabilization*.

As previously mentioned, we entirely restrict ourselves to steady problems and so neglect all temporal strategies. This is due in part to the scarcity of literature on DPG methods for transient problems. However, in this regard, the interested reader may wish to consult [27–30].

*1.3. Camellia*

In our study, we rely heavily upon *Camellia* [3, 4], a C++ toolbox developed by Nathan V. Roberts which uses Sandia's Trilinos library of packages [31]. The general approach in this work follows a trajectory for the DPG method which has been established in the previously mentioned Stokes and Navier-Stokes studies.

Although it is not appropriate to discuss all of the core features here, Camellia is a publicly available software[3] with many tools for rapid implementation of several different finite element methods including discontinuous Galerkin, discontinuous Petrov-Galerkin, hybridizable discontinuous Galerkin [32], and first-order system least-squares [33].

---

[2]However, the linear systems which are obtained with stabilized methods are generally not symmetric and their discrete stability constants can be far from those inherent to the original infinite-dimensional problem.

[3]Available at https://bitbucket.org/nateroberts/camellia.git.

Its mechanisms allow the user to simply provide the bilinear form, boundary conditions, polynomial order, test norm, and material and load data before solving a problem. It also supports shape functions coming from the full exact sequence for *most* standard one-, two-, and three-dimensional elements. Specifically, it supports the conforming and non-conforming two-dimensional quadrilateral $H^1$ and $H(\text{div})$ shape functions, which were used in our work.

### 1.4. Notations and conventions

Throughout this paper, we regularly resort to abstract linear and bilinear operator notation. For these operators, several finite element colloquialisms will be used. For instance, for an abstract bilinear form on Hilbert spaces; $b : U \times V \to \mathbb{R}$, $U$ will denote the *trial space* and $V$ the *test space*.

For some *load*, or continuous linear form, $\ell \in V'$, we are generally interested in the equivalent abstract problems

$$\begin{cases} \text{Find } \mathfrak{u} \in U, \\ b(\mathfrak{u}, \mathfrak{v}) = \ell(\mathfrak{v}), \quad \forall \mathfrak{v} \in V, \end{cases} \iff \begin{cases} \text{Find } \mathfrak{u} \in U, \\ B\mathfrak{u} = \ell, \end{cases} \quad (1.4)$$

where $\langle B\mathfrak{u}, \mathfrak{v} \rangle_{V' \times V} = b(\mathfrak{u}, \mathfrak{v})$ for all $\mathfrak{v} \in V$. The operators $b$ and $B$ are entirely interchangeable and we will often pass between them for simplicity of exposition.

In order to define the specific spaces that will take the places of $U$ and $V$ invoked above, we must define some typical Hilbert spaces. To begin, we define the $L^2$ inner product in a domain $K \subseteq \Omega$ as

$$(u, v)_K = \int_K \text{Tr}(u^\mathsf{T} v) \, dK,$$

where Tr is the usual algebraic trace of a matrix. That is, depending upon whether $u$ and $v$ take scalar, vector, or matrix values, $\text{Tr}(u^\mathsf{T} v)$ will be $uv$, $u \cdot v$, or $u : v$, respectively. The specific Lebesgue spaces that we will need are

$$L^2(K) = \{p : K \to \mathbb{R} \mid \|p\|^2_{L^2(K)} = (p, p)_K < \infty\},$$
$$\boldsymbol{L}^2(K) = \{\mathbf{u} : K \to \mathbb{R}^2 \mid \|\mathbf{u}\|^2_{\boldsymbol{L}^2(K)} = (\mathbf{u}, \mathbf{u})_K < \infty\},$$
$$\boldsymbol{L}^2(K; \mathbb{U}) = \{\boldsymbol{\sigma} : K \to \mathbb{U} \mid \|\boldsymbol{\sigma}\|^2_{\boldsymbol{L}^2(K;\mathbb{U})} = (\boldsymbol{\sigma}, \boldsymbol{\sigma})_K < \infty\},$$

where $\mathbb{U}$ is a subspace of $\mathbb{M}$, the space of $2 \times 2$ real-valued matrices. More explicitly, $\mathbb{U}$ will be allowed to be the symmetric matrices, $\mathbb{S}$, or $\mathbb{M}$ itself. We will use these Lebesgue spaces to construct the solution of the ultraweak formulation of the viscoelastic fluid equations; since no derivatives will be applied to the solution variables, no regularity conditions need to be assumed.

Letting the mesh be denoted $\mathcal{T}$, the norm on the space of discontinuous test velocities is defined

$$\|\mathbf{v}\|^2_{\boldsymbol{H}^1(\mathcal{T})} = \sum_{K \in \mathcal{T}} \left( \|\mathbf{v}\|^2_{\boldsymbol{L}^2(K)} + \|\boldsymbol{\nabla} \mathbf{v}\|^2_{\boldsymbol{L}^2(K;\mathbb{M})} \right),$$

where $\boldsymbol{\nabla}$ denotes the row-wise distributional gradient. Each term in the parentheses above is the square of the standard $H^1$-norm of $\mathbf{v}$, $\|\mathbf{v}\|^2_{\boldsymbol{H}^1(K)}$. The associated Hilbert space of discontinuous test velocities is simply

$$\boldsymbol{H}^1(\mathcal{T}) = \{\mathbf{v} : \Omega \to \mathbb{R}^2 \mid \|\mathbf{v}\|_{\boldsymbol{H}^1(\mathcal{T})} < \infty\}, \qquad (1.5)$$

however, we emphasize to the reader that $\boldsymbol{H}^1(\Omega) \subsetneq \boldsymbol{H}^1(\mathcal{T})$ since only the latter contains functions with jump discontinuities.

All other discontinuous test function spaces are defined similarly.

Let us now consider a *creeping* flow, wherein the conservation of momentum dictates that the Cauchy stress obey the relationship

$$-\boldsymbol{\nabla} \cdot \boldsymbol{\sigma} = \rho \mathbf{f},$$

where $\rho$ is the *mass density*, $\mathbf{f}$ is a *body force density*, and $\boldsymbol{\nabla} \cdot$ denotes the row-wise distributional divergence. Multiplying this equation by a test velocity, $\mathbf{v}$, integrating over a single element, $K$, and then integrating by parts, we obtain

$$(\boldsymbol{\sigma}, \boldsymbol{\nabla} \mathbf{v})_K - \langle \boldsymbol{\sigma} \cdot \hat{\mathbf{n}}, \mathbf{v} \rangle_{\partial K} = (\rho \mathbf{f}, \mathbf{v})_K.$$

Here, we understand $\langle \cdot, \cdot \rangle_{\partial K}$ to indicate (assuming sufficiently smooth variables[4]) an integral over the element boundary.

The DPG method will require that we disassociate the stress variable inside the element with its action on the boundary. To do this, we introduce a new unknown flux variable, $\hat{\mathbf{t}} = \begin{pmatrix} \hat{t}_1 \\ \hat{t}_2 \end{pmatrix}$ —which we call the *traction*—to replace $\boldsymbol{\sigma} \cdot \hat{\mathbf{n}}$. Summing over each element in the mesh, we arrive at the equation

$$(\boldsymbol{\sigma}, \boldsymbol{\nabla}_h \mathbf{v})_\Omega - \langle \hat{\mathbf{t}}, \mathbf{v} \rangle_h = (\rho \mathbf{f}, \mathbf{v})_\Omega.$$

Here, $\boldsymbol{\nabla}_h$ indicates that the gradient is intended element-wise and $\langle \cdot, \cdot \rangle_h$ indicates the accumulation of all related boundary terms.

In the DPG method, trial variables with a circumflex are called *interface variables*. They only have values over the skeleton of the mesh; they are not defined on element interiors. The trial variables defined local to the elements are called *field variables* and in the ultraweak setting are discontinuous and are not well-defined on the mesh skeleton. For a further example, we will eventually use $\mathbf{u}$ to represent the *field* velocity and $\hat{\mathbf{u}}$ to represent the *interface* velocity.

Interface variables in the DPG method act as Lagrange multipliers in response to the discontinuity of the test space. In an ultraweak formulation, the interface variables are the only variables which have continuity across element boundaries. That is, we allow discontinuity in the field variables in an ultraweak formulation.

We call interface variables which replace terms involving the outward facing normal *flux* variables (e.g. $\hat{\mathbf{t}}$), and interface variables simply representing restriction to the boundary *trace* variables (e.g. $\hat{\mathbf{u}}$). This distinction is important because in a 2D computer implementation, flux variable continuity is enforced only at edge nodes while trace variable continuity is also enforced at vertex nodes. The reasons for this are theoretical and discussed in [3].

---
[4]Lower regularity assumptions require this operation to be identified with a duality pairing. See [34].

*1.5. Outline*

The next section gives a compact outline of the features of this non-standard finite element method most important for our study on viscoelastic fluids. Then, once the essential components of the DPG method have been established, in Section 3 we provide a precise ultraweak variational formulation that we will use in the benchmark analysis. In Section 4 we verify our implementation by comparing with the existing literature. The paper closes in Section 5 with concluding remarks and a brief discussion of aspirations for future work.

## 2. DPG: the main ideas

Abstractly, the method is derived as follows. Let $b : U \times V \to \mathbb{R}$ be a bilinear form on Hilbert spaces, $U$ and $V$. After a discretization of the trial space has been chosen—or equivalently, after a set of computable solutions, $U_h \subset U$, has been fixed—we seek the optimal discrete solution, $\mathfrak{u}_h^{\text{opt}}$, of the residual minimization problem

$$\mathfrak{u}_h^{\text{opt}} = \arg\min_{\mathfrak{u}_h \in U_h} \|B\mathfrak{u}_h - \ell\|_{V'}^2. \quad (2.1)$$

This residual optimality translates into optimality of the error in the induced energy norm, $\|\mathfrak{u}\|_{U^{\text{ind.}}} = \|B\mathfrak{u}\|_{V'}$, since the (induced) energy error can be written as

$$\|\mathfrak{u} - \mathfrak{u}_h\|_{U^{\text{ind.}}} = \|B(\mathfrak{u} - \mathfrak{u}_h)\|_{V'} = \|B\mathfrak{u}_h - \ell\|_{V'}.$$

Such a derived norm is determined only by the operator, $B$, and the test norm being used, $\|\cdot\|_V$.

We can easily characterize the optimal solution as the unique stationary solution of (2.1). First, however, it is helpful to rewrite

$$\|B\mathfrak{u}_h - \ell\|_{V'}^2 = \langle B\mathfrak{u}_h - \ell, \mathcal{R}_{V'}(B\mathfrak{u}_h - \ell)\rangle, \quad (2.2)$$

where $\mathcal{R}_V : V \to V'$ is the Riesz operator defined by

$$\langle \mathcal{R}_V \mathfrak{v}, \delta\mathfrak{v}\rangle = (\mathfrak{v}, \delta\mathfrak{v})_V \quad \forall \delta\mathfrak{v} \in V, \quad (2.3)$$

and $(\cdot, \cdot)_V$ is the inner product on $V$. The Riesz map is an isometry—that is, $\|\mathcal{R}_V v\|_{V'} = \|v\|_V$—so that $\mathcal{R}_{V'} = \mathcal{R}_V^{-1}$. Therefore, invoking (2.2), the optimal solution must vanish in the first variation of (2.1),

$$\langle B\mathfrak{u}_h^{\text{opt}} - \ell, \mathcal{R}_V^{-1} B\delta\mathfrak{u}_h\rangle = 0 \quad \forall \delta\mathfrak{u}_h \in U_h,$$

or, equivalently,

$$b(\mathfrak{u}_h^{\text{opt}}, \mathcal{R}_V^{-1} B\delta\mathfrak{u}_h) = \ell(\mathcal{R}_V^{-1} B\delta\mathfrak{u}_h)) \quad \forall \delta\mathfrak{u}_h \in U_h. \quad (2.4)$$

If we define the DPG inner product as $a(\mathfrak{u}, \delta\mathfrak{u}) = \langle B\mathfrak{u}, \mathcal{R}_V^{-1} B\delta\mathfrak{u}\rangle$ for all $\mathfrak{u}, \delta\mathfrak{u} \in U$, then symmetry of $a$ is obvious:

$$a(\mathfrak{u}, \delta\mathfrak{u}) = \langle B\mathfrak{u}_h, \mathcal{R}_V^{-1} B\delta\mathfrak{u}_h\rangle = \langle \mathcal{R}_V^{-1} B\mathfrak{u}_h, B\delta\mathfrak{u}_h\rangle = a(\delta\mathfrak{u}, \mathfrak{u}).$$

Notice that computing the optimal solution involves evaluating the inverse of the Riesz map. Unfortunately, computing the inverse of a global operator is prohibitively expensive in general scenarios. However, this can be made practical when using a test space which is discontinuous across element boundaries. Such test spaces can be decomposed into a finite orthogonal direct sum, $V^{\text{DPG}} = \bigoplus_{K \in \mathcal{T}} V_K$. Here, it is helpful to envision each $V_K$ as an element-wise test space, where $K$ is a given element in the mesh, $\mathcal{T}$. For instance, consider the space of test velocities in (1.5) where $V^{\text{DPG}} = \boldsymbol{H}^1(\mathcal{T})$ and each $V_K = \boldsymbol{H}^1(K)$. With such a decomposition, we see that the Riesz map can be defined through its action on individual elements. This is because for every element-local test function, $\mathfrak{v}_K|_K \in V_K$ and $\mathfrak{v}_K|_{\Omega \setminus K} = 0$, and all possible $\delta\mathfrak{v}^{\text{DPG}} = \sum_{\widetilde{K} \in \mathcal{T}} \delta\mathfrak{v}_{\widetilde{K}}$,

$$\langle \mathcal{R}_{V^{\text{DPG}}} \mathfrak{v}_K, \delta\mathfrak{v}^{\text{DPG}}\rangle = \sum_{\widetilde{K} \in \mathcal{T}} (\mathfrak{v}_K, \delta\mathfrak{v}_{\widetilde{K}})_{V^{\text{DPG}}}$$
$$= (\mathfrak{v}_K, \delta\mathfrak{v}_K)_{V^{\text{DPG}}} = \langle \mathcal{R}_{V_K} \mathfrak{v}_K, \delta\mathfrak{v}_K\rangle.$$

This leads to locality of the operator and its inverse, $\mathcal{R}_{V^{\text{DPG}}}^{-1} = \bigoplus_{K \in \mathcal{T}} \mathcal{R}_{V_K}^{-1}$, which can therefore be efficiently approximated on-the-fly and *in parallel*. With this locality, we can also decompose the residual into a *single* sum over the elements of the mesh,

$$\|B\mathfrak{u}_h - \ell\|_{V^{\text{DPG}'}}^2 = \sum_{K \in \mathcal{T}} \langle B\mathfrak{u}_h - \ell, \mathcal{R}_{V_K}^{-1}(B\mathfrak{u}_h - \ell)\rangle$$
$$= \sum_{K \in \mathcal{T}} \|B\mathfrak{u}_h - \ell\|_{V_K'}^2 = \sum_{K \in \mathcal{T}} \eta_{V_K}^2. \quad (2.5)$$

Notable here is that each individual term on the right hand side of (2.5), induces an *a posterior* error indicator which we denote for an element $K$ as $\eta_{V_K}$. We will use such estimates of local error contributions to infer an intrinsic refinement strategy which we call the *energy strategy*.

As mentioned in Section 1.4, due to the discontinuous test spaces, a DPG solution usually involves interface variables. By decomposing a DPG trial variable into two terms, $\mathfrak{u}^{\text{DPG}} = (\mathfrak{u}^{\text{fld.}}, \hat{\mathfrak{u}})$, we can also decompose the DPG bilinear form, $b^{\text{DPG}}$, into two separate contributions,

$$b^{\text{DPG}}((\mathfrak{u}^{\text{fld.}}, \hat{\mathfrak{u}}), \mathfrak{v}^{\text{DPG}}) = b^{\text{fld.}}(\mathfrak{u}^{\text{fld.}}, \mathfrak{v}^{\text{DPG}}) + \hat{b}(\hat{\mathfrak{u}}, \mathfrak{v}^{\text{DPG}}).$$

Here, $\mathfrak{u}^{\text{fld.}}$ is solution field variable with values across the entire domain, $\hat{\mathfrak{u}}$ is an interface variable with values only upon the mesh skeleton, and $\mathfrak{v}^{\text{DPG}} = \sum_{K \in \mathcal{T}} \mathfrak{v}_K \in V^{\text{DPG}}$ is a test function which is allowed to be discontinuous across all element boundaries.

*2.1. Stability*

Stability in a finite element method essentially determines—in a quantifiable way—whether the discrete problem is well-posed under compatible loads [35]. For arbitrary discrete trial and test spaces, $U_h \subset U$ and $V_h \subset V$, the amount of stability in a finite element method is measured by the so-called *discrete inf-sup constant*,

$$\gamma_h = \inf_{\mathfrak{u}_h \in U_h} \sup_{\mathfrak{v}_h \in V_h} \frac{|b(\mathfrak{u}_h, \mathfrak{v}_h)|}{\|\mathfrak{u}_h\|_U \|\mathfrak{v}_h\|_V}.$$

The burden here is not only that without stability—that is, if $\gamma_h = 0$—the method *will not converge*, but that the smaller the value of $\gamma_h$, the larger the difference can be between the discrete and the exact solution.

The generation and analysis of stable finite element methods can be very challenging; thus a wide array of techniques have been developed for such pursuits. Although many such techniques have become very successful, a large amount of time—that is, time before code is ever run successfully—can still be devoted to developing stabilizations of discretizations each time a new problem is to be solved. The DPG method is thus very appealing in that it is an intrinsically stable finite element method. i.e. *No stabilization techniques ever need to be applied.* The reason for this is actually quite simple: by allowing for a given $\mathfrak{u}_h \in U_h$, $\mathfrak{v}_h = \mathcal{R}_V^{-1} B \mathfrak{u}_h \in V$,

$$\frac{b(\mathfrak{u}_h, \mathfrak{v}_h)}{\|\mathfrak{v}_h\|_V} = \frac{(\mathcal{R}_V^{-1} B \mathfrak{u}_h, \mathcal{R}_V^{-1} B \mathfrak{u}_h)_V}{\|\mathcal{R}_V^{-1} B \mathfrak{u}_h\|_V}$$
$$= \|\mathcal{R}_V^{-1} B \mathfrak{u}_h\|_V = \|B \mathfrak{u}_h\|_{V'} = \sup_{\mathfrak{v} \in V} \frac{|b(\mathfrak{u}_h, \mathfrak{v})|}{\|\mathfrak{v}\|_V}.$$

Therefore, recalling (2.4), we immediately establish that for $V_h = \mathcal{R}_V^{-1} B(U_h)$,

$$\gamma_h = \inf_{\mathfrak{u}_h \in U_h} \sup_{\mathfrak{v}_h \in V_h} \frac{|b(\mathfrak{u}_h, \mathfrak{v}_h)|}{\|\mathfrak{u}_h\|_U \|\mathfrak{v}_h\|_V} = \inf_{\mathfrak{u}_h \in U_h} \sup_{\mathfrak{v} \in V} \frac{|b(\mathfrak{u}_h, \mathfrak{v})|}{\|\mathfrak{u}_h\|_U \|\mathfrak{v}\|_V}.$$

In fact, the right-most expression bounds the *continuous inf-sup constant*, $\gamma$.[5] Thus, if $\gamma$ is bounded away from zero, then $\gamma_h$ must be as well. Notably, the requirement that $\gamma > 0$ is very modest: this will be true for any well-posed problem so it must hold if a unique solution exists.

Practical algorithms for leveraging this observation have eluded researchers due to the computational complexity of inverting the Riesz map. Of course, we still do not exactly invert what is still an *infinite-dimensional* operator, $\mathcal{R}_{V^{\text{DPG}}}$. Instead, by approximating its explicit inverse to a tunable accuracy at the element level, we understand that we will eventually reach a point where stability is achieved. Generally, this threshold is low [36, 37] and in our computational experiments we do indeed observe stability quickly and reliably with this approach.

### 2.2. Optimal test norms

Thus far, we have left the norm on $V$ be user-defined. However, some choices must be more desirable than others, especially recalling that the choice of test norm affects which solution we arrive with via (2.1). We should also keep in mind that the error indicator $\eta_{V_K}$ is determined by the test norm we use, and therefore the choice of test norm will also affect the refinement patterns that we produce.

Inspiration comes from attempting to construct a test norm such that the induced energy norm, $\|\cdot\|_{U^{\text{ind.}}} = \|B \cdot \|_{V'}$, coincides with a particular norm of interest, $\|\cdot\|_U$. We define the *optimal test norm* of $\mathfrak{v} \in V$ [38] to be

$$\|\mathfrak{v}\|_{V^{\text{opt}}} = \sup_{\mathfrak{u} \in U} \frac{|\langle B \mathfrak{u}, \mathfrak{v} \rangle|}{\|\mathfrak{u}\|_U} = \|B' \mathfrak{v}\|_{U'}. \quad (2.6)$$

---
[5]Explicitly, this bound is

$$\gamma_h = \inf_{\mathfrak{u}_h \in U_h} \sup_{\mathfrak{v} \in V} \frac{|b(\mathfrak{u}_h, \mathfrak{v})|}{\|\mathfrak{u}_h\|_U \|\mathfrak{v}\|_V} \geq \inf_{\mathfrak{u} \in U} \sup_{\mathfrak{v} \in V} \frac{|b(\mathfrak{u}, \mathfrak{v})|}{\|\mathfrak{u}\|_U \|\mathfrak{v}\|_V} = \gamma.$$

Likewise, $\mathcal{R}_{V^{\text{opt}}} = B' \mathcal{R}_U^{-1} B$. Observe that with this definition we can write

$$\|\mathfrak{u}\|_U = \sup_{\mathfrak{u}' \in U'} \frac{|\langle \mathfrak{u}, \mathfrak{u}' \rangle|}{\|\mathfrak{u}'\|_{U'}} = \sup_{\mathfrak{v} \in V} \frac{|\langle \mathfrak{u}, B' \mathfrak{v} \rangle|}{\|B' \mathfrak{v}\|_{U'}}$$
$$= \sup_{\mathfrak{v} \in V} \frac{|\langle B \mathfrak{u}, \mathfrak{v} \rangle|}{\|\mathfrak{v}\|_{V^{\text{opt}}}} = \|B \mathfrak{u}\|_{(V^{\text{opt}})'},$$

for all $\mathfrak{u} \in U$. Here, we recall that since $B$ is injective[6] then $B' : V \to U'$ is surjective and we arrive with the second equality above. Therefore, by using the corresponding optimal test norm, the solution error in our *chosen trial norm* can always be written

$$\|\mathfrak{u} - \mathfrak{u}_h\|_U = \|B \mathfrak{u}_h - \ell\|_{(V^{\text{opt}})'}. \quad (2.7)$$

In the ultraweak setting, $\mathcal{R}_{V^{\text{opt}}}$ is often nearly fully computable. Indeed, $b^{\text{fld.}}(\mathfrak{u}^{\text{fld.}}, \mathfrak{v}^{\text{DPG}}) = (\mathfrak{u}^{\text{fld.}}, \mathcal{B}^* \mathfrak{v}^{\text{DPG}})_{L^2(\Omega)}$ for a unique operator $\mathcal{B}^* : V \to L^2$, so if the norm of interest comes from an $L^2$ norm on the field variables, $\|\mathfrak{u}^{\text{fld.}}\|_{L^2}$, then

$$\|\mathfrak{v}\|_{V^{\text{opt,fld.}}} = \sup_{\mathfrak{u}^{\text{fld.}} \in U^{\text{fld.}}} \frac{|b^{\text{fld.}}(\mathfrak{u}^{\text{fld.}}, \mathfrak{v})|}{\|\mathfrak{u}^{\text{fld.}}\|_{U^{\text{fld.}}}} = \|\mathcal{B}^* \mathfrak{v}\|_{L^2(\Omega)}.$$

In the definition of this global (semi-)norm above, we have ignored the contributions of the interface terms which would make it an optimal test norm on $V^{\text{DPG}}$. Therefore, it is of little surprise that the related restriction to each element, $\|\mathcal{B}^* \mathfrak{v}^{\text{DPG}}\|_{L^2(K)}$, is not a norm and furthermore, $\|\mathfrak{v}\|_{V^{\text{opt,fld.}}}$ is only a semi-norm on the discontinuous test space, $V^{\text{DPG}}$.

We rectify this issue incorporating $L^2$ terms to construct what we call the *(adjoint) graph norm*,

$$\|\mathfrak{v}^{\text{DPG}}\|_{V^{\text{graph}}}^2 = \sum_{K \in \mathcal{T}} \left( \|\mathcal{B}^* \mathfrak{v}^{\text{DPG}}\|_{L^2(K)}^2 + \|\mathfrak{v}^{\text{DPG}}\|_{L^2(K)}^2 \right), \quad (2.8)$$

which can be easily discretized. In several fluid flow problems, the graph norm has proved effective [24, 25] and we will also use it here.

### 2.3. Nonlinear Problems

Our method for applying DPG to nonlinear problems follows a well-established trajectory explored in many papers [39]. Most notably, steady Stokes flow problems were first handled with the DPG method in [24] and using insight from that problem, a nonlinear algorithm was constructed to handle steady Navier-Stokes problems in [25]. Our approach here builds directly upon this previous work.

The general strategy is constructed from the standard Gauss-Newton algorithm. Our approach is to apply the DPG methodology to the linearized problem and successively minimize the residual of the linearized problem at each iteration. Specifically, we begin with an initial guess for the solution, $\mathfrak{u}_0$, and about this point we construct a linearized problem for a solution increment, $\Delta \mathfrak{u}$, *viz.*,

$$\begin{cases} \text{Find } \Delta \mathfrak{u} \in U, \\ b_{\text{lin.}}[\mathfrak{u}_0](\Delta \mathfrak{u}, \mathfrak{v}) = \ell_{\text{lin.}}[\mathfrak{u}_0](\mathfrak{v}), \quad \forall \mathfrak{v} \in V. \end{cases}$$

---
[6]This follows if (1.4) is well-posed.

We define the linearized form to be the derivative of the nonlinear form, $b_{\text{lin.}}[\mathfrak{u}_0](\Delta\mathfrak{u},\mathfrak{v}) = D_\mathfrak{u} b_{\text{nl.}}[\mathfrak{u}_0](\Delta\mathfrak{u},\mathfrak{v})$ and define the forcing term to be the (nonlinear) residual $\ell_{\text{lin.}}[\mathfrak{u}_0](\mathfrak{v}) = \ell(\mathfrak{v}) - b_{\text{nl.}}(\mathfrak{u}_0,\mathfrak{v})$. During each iteration, we increment the solution in the optimal direction $\mathfrak{u}_0 \mapsto \mathfrak{u}_0 + \Delta\mathfrak{u}_h^{\text{opt}}$ where $\Delta\mathfrak{u}_h^{\text{opt}}$. Here, $\Delta\mathfrak{u}_h^{\text{opt}}$ minimizes the residual over the trial space, $U_h \subset U$,

$$\Delta\mathfrak{u}_h^{\text{opt}} = \underset{\Delta\mathfrak{u}_h \in U_h}{\arg\min} \left\| B_{\text{lin.}}[\mathfrak{u}_0]\Delta\mathfrak{u}_h - \ell_{\text{lin.}}[\mathfrak{u}_0] \right\|_{V'},$$

and $B_{\text{lin.}}[\mathfrak{u}_0]$ is defined through $b_{\text{lin.}}[\mathfrak{u}_0]$ as in (1.4). An important feature of our work is that, when using the graph norm, the test space will inherit a norm which is updated *along with* the solution increment, $\mathfrak{u}_0$. Updating the test norm, in this way allows us to ensure discrete stability of each linearized problem.

### 2.4. Implementing DPG

The fine details of the numerical implementation of the DPG method are described in [3, 24, 40] but, ultimately, the inverse of the Riesz map on each element must be approximated. Therefore, all of the previous results concerning stability and optimality stay true only in an asymptotic sense.

The method is usually implemented using an order-elevated polynomial space in $V^{\text{enr}} = \bigoplus_{K \in \mathcal{T}} V_K^{\text{enr}} \subset V^{\text{DPG}}$ where, at each element, the polynomial order has been increased above that of the trial space, $p$, by a small increment $\mathrm{d}p$. Numerical evidence supports that $\mathrm{d}p$ can be very small in practice—sometimes, even in 3D, $\mathrm{d}p = 1$ gives reliable results. Theoretical results, however, suggest that the increment increases with the spatial dimension [41].

We define $\mathcal{P}^p(K)$ as the space of uniform $p$-order polynomials on the subset of the domain $K \subset \Omega$. In an ultraweak discretization, we draw each component of the field variables from a discontinuous polynomial space of uniform order, $p$, denoted $\mathcal{P}^p(\mathcal{T}) = \{f|_K \in \mathcal{P}^p(K) \mid K \in \mathcal{T}\}$. The interface variables are drawn from similarly-defined spaces, but must be continuous across each edge. For theoretical reasons, we choose the trace variables to be in $\mathcal{P}^{p+1}(E)$ and the flux variables from $\mathcal{P}^p(E)$ for each edge, $E$, in the mesh. In other implementations, the interface variables are drawn from restrictions of standard exact sequence conforming polynomial spaces such as those developed in [42, 43]; however, this has not been implemented here.

For the test functions, we construct the discretization from non-conforming analogues of the standard $H^1$- or $H(\text{div})$-conforming polynomial spaces [42], depending, respectively, upon the functional setting of each component of the test function. In each case, we use a polynomial space coming from the exact sequence of order $p + \mathrm{d}p$.

We emphasize to the reader that both the Riesz operator as well as the bilinear form need to be discretized with this method. Therefore, because the Riesz operator also needs to be inverted at each element, an LU-decomposition and a back-substitution must be performed, followed by matrix-matrix multiplication.

For a more dedicated account of the practical implementation of the method, especially with Camellia, we refer the interested reader to [3, 4].

## 3. Formulation

In this section, we formally derive the variational formulation we will use for the Oldroyd-B model with Navier-Stokes coupling and then the corresponding test space adjoint graph norm. The derivations for the Giesekus model, $\alpha > 0$, and the Oldroyd-B model with only Stokes coupling are left to the reader.

### 3.1. A mesh-dependent ultraweak formulation

By introducing the velocity gradient, $\mathbf{L} = \boldsymbol{\nabla}\mathbf{u}$, as an independent variable, we can write the entire coupled steady system as

$$\rho(\mathbf{L}\mathbf{u}) + \nabla p - \eta_{\text{S}}\boldsymbol{\nabla}\cdot\mathbf{L} - \boldsymbol{\nabla}\cdot\mathbf{T} = \rho\mathbf{f}, \quad (3.1)$$

$$\mathbf{L} - \boldsymbol{\nabla}\mathbf{u} = \mathbf{0}, \quad (3.2)$$

$$\boldsymbol{\nabla}\cdot\mathbf{u} = 0, \quad (3.3)$$

$$\mathbf{T} + \lambda\mathcal{L}_u\mathbf{T} - \eta_{\text{P}}(\mathbf{L} + \mathbf{L}^\mathsf{T}) = \mathbf{0}, \quad (3.4)$$

where we have introduced the *autonomous Lie derivative* [44, Section 1.6] of $\mathbf{T}$ in the direction of $\mathbf{u}$, $\mathcal{L}_\mathbf{u}\mathbf{T} = \mathfrak{L}_\mathbf{u}\mathbf{T} - \frac{\partial\mathbf{T}}{\partial t}$. Taking into account (3.2) and (3.3), this can be expressed as

$$\mathcal{L}_u\mathbf{T} = \boldsymbol{\nabla}\cdot(\mathbf{T}\otimes\mathbf{u}) - \mathbf{L}\,\mathbf{T} - \mathbf{T}\,\mathbf{L}^\mathsf{T}.$$

The solution components $\mathbf{u}$, $p$, $\mathbf{L}$, and $\mathbf{T}$ will constitute all of the field variables in our ultraweak formulation. Formally testing (3.1) with smooth vector fields, $\mathbf{v}$, (3.2) with smooth 2-tensors, $\mathbf{M}$, (3.3) with smooth functions,[7] $q$, and (3.4) with smooth *symmetric* tensors, $\mathbf{S}$, we see that after integration by parts over each element and summing each equation together over each element of the mesh, we arrive at the nonlinear form

$$\begin{aligned}b_{\text{nl.}}^{\text{fld.}}&((\mathbf{u},p,\mathbf{L},\mathbf{T}),(\mathbf{v},q,\mathbf{M},\mathbf{S}))\\&= \rho(\mathbf{L}\mathbf{u},\mathbf{v}) - (p,\nabla_h\cdot\mathbf{v}) + \eta_{\text{S}}(\mathbf{L},\boldsymbol{\nabla}_h\mathbf{v}) + (\mathbf{T},\boldsymbol{\nabla}_h\mathbf{v})\\&\quad + (\mathbf{L},\mathbf{M}) + (\mathbf{u},\boldsymbol{\nabla}_h\cdot\mathbf{M}) - (\mathbf{u},\nabla_h q)\\&\quad + (\mathbf{T},\mathbf{S}) - \lambda(\mathbf{T}\otimes\mathbf{u},\boldsymbol{\nabla}_h\mathbf{S}) - 2\lambda(\mathbf{L}\,\mathbf{T},\mathbf{S}) - 2\eta_{\text{P}}(\mathbf{L},\mathbf{S}).\end{aligned}$$

Here, for brevity, we have discarded the subscript "$\Omega$" from the $L^2$ inner product notation. In order to arrive at an equation without boundary terms, we have assumed the test variables to vanish on the domain boundary, $\partial\Omega$. If, however, the test functions are unrestricted on the boundary of the domain and are allowed to be discontinuous across element interfaces then the additional contributions collect together to form the interface terms

$$\begin{aligned}\hat{b}&((\hat{\mathbf{u}},\hat{\mathbf{t}},\hat{\mathbf{j}}),(\mathbf{v},q,\mathbf{M},\mathbf{S}))\\&= -\langle\hat{\mathbf{t}},\mathbf{v}\rangle_h - \langle\hat{\mathbf{u}},\mathbf{M}\hat{\mathbf{n}}\rangle_h + \langle\hat{\mathbf{u}}\cdot\hat{\mathbf{n}},q\rangle_h + \lambda\langle\hat{\mathbf{j}},\mathbf{S}\rangle_h.\end{aligned}$$

If we assume a smooth solution near the boundary of each element of the mesh, $K \in \mathcal{T}$, the new interface solution variables can be identified with the restriction of the velocity field, $\hat{\mathbf{u}} = \mathbf{u}|_{\partial K}$, the normal flux of the Cauchy stress, $\hat{\mathbf{t}} = \boldsymbol{\sigma}|_{\partial K}\hat{\mathbf{n}}_K$,

---

[7]Note that we could instead enforce $\text{Tr}(\mathbf{L}) = 0$.

and the normal flux in the hybrid variable $\mathbf{T} \otimes \mathbf{u}$, $\hat{\mathbf{j}} = (\mathbf{u}|_{\partial K} \cdot \hat{\mathbf{n}}_K)\mathbf{T}|_{\partial K}$.[8] In this case, the entire nonlinear form becomes

$$b_{\text{nl.}}((\mathbf{u}, p, \mathbf{L}, \mathbf{T}, \hat{\mathbf{u}}, \hat{\mathbf{t}}, \hat{\mathbf{j}}), (\mathbf{v}, q, \mathbf{M}, \mathbf{S}))$$
$$= b_{\text{nl.}}^{\text{fld.}}((\mathbf{u}, p, \mathbf{L}, \mathbf{T}), (\mathbf{v}, q, \mathbf{M}, \mathbf{S})) + \hat{b}((\hat{\mathbf{u}}, \hat{\mathbf{t}}, \hat{\mathbf{j}}), (\mathbf{v}, q, \mathbf{M}, \mathbf{S})).$$

We note that $\hat{b}(\cdot, \cdot)$ is in fact linear in each of its arguments, however, $b_{\text{nl.}}^{\text{fld.}}(\cdot, \cdot)$ is only linear in the pressure variable, $p$.

Following the strategy of Section 2.3, but from now on abandoning the prefix $\Delta$ from each solution variable as well as the subscript-$h$ denoting element-wise differentiation, we find that

$$b_{\text{lin.}}^{\text{fld.}}[\mathbf{u}_0, \mathbf{L}_0, \mathbf{T}_0]((\mathbf{u}, p, \mathbf{L}, \mathbf{T}), (\mathbf{v}, q, \mathbf{M}, \mathbf{S}))$$
$$= D_{\mathbf{u}} b_{\text{nl.}}^{\text{fld.}}[\mathbf{u}_0, 0, \mathbf{L}_0, \mathbf{T}_0]((\mathbf{u}, p, \mathbf{L}, \mathbf{T}), (\mathbf{v}, q, \mathbf{M}, \mathbf{S}))$$
$$= \rho(\mathbf{L}_0 \mathbf{u} + \mathbf{L} \mathbf{u}_0, \mathbf{v}) - (p, \nabla \cdot \mathbf{v}) + \eta_{\text{S}}(\mathbf{L}, \nabla \mathbf{v})$$
$$+ (\mathbf{T}, \nabla \mathbf{v}) + (\mathbf{L}, \mathbf{M}) + (\mathbf{u}, \nabla \cdot \mathbf{M}) - (\mathbf{u}, \nabla q)$$
$$+ (\mathbf{T}, \mathbf{S}) - \lambda(\mathbf{T}_0 \otimes \mathbf{u} + \mathbf{T} \otimes \mathbf{u}_0, \nabla \mathbf{S})$$
$$- 2\lambda(\mathbf{L}_0 \mathbf{T} + \mathbf{L} \mathbf{T}_0, \mathbf{S}) - 2\eta_{\text{P}}(\mathbf{L}, \mathbf{S})$$
$$= (\mathbf{u}, \rho \mathbf{L}_0^\top \mathbf{v} - \nabla q + \nabla \cdot \mathbf{M} - \lambda \nabla \mathbf{S} : \mathbf{T}_0)$$
$$+ (\mathbf{L}, \eta_{\text{S}} \nabla \mathbf{v} + \rho \mathbf{v} \otimes \mathbf{u}_0 + \mathbf{M} - 2\eta_{\text{P}} \mathbf{S} - 2\lambda \mathbf{S} \mathbf{T}_0)$$
$$- (p, \nabla \cdot \mathbf{v}) + (\mathbf{T}, \nabla \mathbf{v} + \mathbf{S} - \lambda(\mathbf{u}_0 \cdot \nabla)\mathbf{S} - 2\lambda \mathbf{L}_0^\top \mathbf{S}),$$

where $\nabla \mathbf{S} : \mathbf{T}_0 = \sum_{i,j,k} (\partial_k \mathbf{S}_{ij})(\mathbf{T}_0)_{ij} \mathbf{e}_k$, and also

$$\ell_{\text{nl.}}[\mathbf{u}_0, \mathbf{L}_0, \mathbf{T}_0](\mathbf{v}, q, \mathbf{M}, \mathbf{S})$$
$$= \ell_{\text{lin.}}(\mathbf{v}, q, \mathbf{M}, \mathbf{S}) - b_{\text{nl.}}((\mathbf{u}_0, 0, \mathbf{L}_0, \mathbf{T}_0, \mathbf{0}, \mathbf{0}, \mathbf{0}), (\mathbf{v}, q, \mathbf{M}, \mathbf{S}))$$
$$= \rho(\mathbf{f}, \mathbf{v}) - b_{\text{nl.}}^{\text{fld.}}((\mathbf{u}_0, 0, \mathbf{L}_0, \mathbf{T}_0), (\mathbf{v}, q, \mathbf{M}, \mathbf{S})),$$

by exploiting the fact that $\hat{b}$ is bilinear.

### 3.2. *Derivation of the graph norm*

Suppressing notation for the domain, $\Omega$, we choose to minimize

$$\|(\mathbf{u}, p, \mathbf{L}, \mathbf{T})\|_U^2$$
$$= \|\eta l_0^{-1} \mathbf{u}\|_{L^2}^2 + \|p\|_{L^2}^2 + \|\eta_{\text{S}} \mathbf{L}\|_{L^2(\mathbb{M})}^2 + \|\mathbf{T}\|_{L^2(\mathbb{S})}^2, \quad (3.5)$$

where $\eta = \eta_{\text{P}} + \eta_{\text{S}}$ is the *total viscosity* and $l_0$ is a length scale to be set for the problem. Notably, in the derivation for the Giesekus and Oldroyd-B models with Stokes flow coupling, we choose a combination of scales such that the weight on the $\|\mathbf{u}\|_{L^2}$-term evaluates to $1 \, \text{kg}^2/(\text{m}^4 \cdot \text{s}^2)$.

Using Lemma A.1, we arrive at an explicit characterization of

$$\|\mathfrak{v}^{\text{DPG}}\|_{V_{\mathbf{u}_0}^{\text{opt,fld.}}}^2 = \left( \sup_{\mathbf{u}^{\text{fld.}} \in U^{\text{fld.}}} \frac{b_{\text{lin.}}^{\text{fld.}}[\mathbf{u}_0](\mathbf{u}^{\text{fld.}}, \mathfrak{v}^{\text{DPG}})}{\|\mathbf{u}\|_U} \right)^2$$

$$= \left( \sup_{\mathbf{u} \in L^2} \frac{(\mathbf{u}, \rho \mathbf{L}_0^\top \mathbf{v} - \nabla q + \nabla \cdot \mathbf{M} - \lambda \nabla \mathbf{S} : \mathbf{T}_0)}{\frac{\eta}{l_0} \|\mathbf{u}\|_{L^2}} \right)^2$$
$$+ \left( \sup_{\mathbf{L} \in L^2(\mathbb{M})} \frac{(\mathbf{L}, \eta_{\text{S}} \nabla \mathbf{v} + \rho \mathbf{v} \otimes \mathbf{u}_0 + \mathbf{M} - 2\eta_{\text{P}} \mathbf{S} - 2\lambda \mathbf{S} \mathbf{T}_0)}{\eta_{\text{S}} \|\mathbf{L}\|_{L^2(\mathbb{M})}} \right)^2$$
$$+ \left( \sup_{p \in L^2} \frac{-(p, \nabla \cdot \mathbf{v})}{\|p\|_{L^2}} \right)^2$$
$$+ \left( \sup_{\mathbf{T} \in L^2(\mathbb{S})} \frac{(\mathbf{T}, \nabla \mathbf{v} + \mathbf{S} - \lambda(\mathbf{u}_0 \cdot \nabla)\mathbf{S} - 2\lambda \mathbf{L}_0^\top \mathbf{S})}{\|\mathbf{T}\|_{L^2(\mathbb{S})}} \right)^2$$
$$= l_0^2 \eta^{-2} \|\rho \mathbf{L}_0^\top \mathbf{v} - \nabla q + \nabla \cdot \mathbf{M} - \lambda \nabla \mathbf{S} : \mathbf{T}_0\|_{L^2}^2$$
$$+ \eta_{\text{S}}^{-2} \|\eta_{\text{S}} \nabla \mathbf{v} + \rho \mathbf{v} \otimes \mathbf{u}_0 + \mathbf{M} - 2\eta_{\text{P}} \mathbf{S} - 2\lambda \mathbf{S} \mathbf{T}_0\|_{L^2(\mathbb{M})}^2$$
$$+ \|\nabla \cdot \mathbf{v}\|_{L^2}^2 + \|\nabla \mathbf{v} + \mathbf{S} - \lambda(\mathbf{u}_0 \cdot \nabla)\mathbf{S} - 2\lambda \mathbf{L}_0^\top \mathbf{S}\|_{L^2(\mathbb{S})}^2,$$

and the graph norm is simply

$$\|\mathfrak{v}^{\text{DPG}}\|_{V_{\mathbf{u}_0}^{\text{graph}}}^2 = \|\mathfrak{v}^{\text{DPG}}\|_{V_{\mathbf{u}_0}^{\text{opt,fld.}}}^2 + \|\mathfrak{v}^{\text{DPG}}\|_{L^2}^2. \quad (3.6)$$

## 4. Confined cylinder benchmark

We implemented the method above using Camellia [3, 4]; our implementation is available as part of Camellia's current stable release. The verification of our method was performed on the so-called confined cylinder problem, for which many results are available in the literature [9–12, 45–49]. This is a simple two-dimensional problem where the fluid passes through a narrow channel with a centrally placed cylinder impeding its flow. The ratio of the cylinder radius to the channel width to the length of the domain is taken to be precisely $1 : 2 : 15$, as depicted in Figure 4.1.

In this problem, the length of the channel is understood to be sufficiently large for outflow conditions to have little effect upon the quantity of interest which is the drag coefficient,

$$\mathfrak{K} = \frac{1}{\eta \bar{u}} \int_{\Gamma_{\text{C}}} (\boldsymbol{\sigma} \hat{\mathbf{n}}) \cdot \mathbf{e}_1. \quad (4.1)$$

Here, $\bar{u}$ is the average inflow velocity and $\Gamma_{\text{C}}$ is the entire cylinder. Notice that we reduced the computational expense by considering only half of the flow domain with reflectively symmetric solutions. We thus compute on the upper-half of $\Omega$, which we denote $\Omega^+$. Because our interest is driven by benchmarking, we do not present a qualitative analysis of our solutions with the common viscometric functions. In this regard, however, our results did closely agree with those reported in [47].

This section continues by further describing the boundary conditions applied to this problem. We then compare our results to the literature for a various values of the Weissenberg number, Wi, various values of the Reynolds number, Re, and various values of the mobility factor, $\alpha$. Of these comparisons, the first dedicated subsection examines the Oldroyd-B model without inertial effects, $\text{Re} = 0$. Here the most extensive results appear in the literature and our comparisons are the most thorough. It is also in this subsection that we analyze our adaptive mesh refinement strategy and investigate the influence of this strategy on our estimates of $\mathfrak{K}$. In Sections 4.3 and 4.4 we use our techniques to examine the effects of non-zero Reynolds numbers followed by non-zero mobility factors. Comparisons are given in both cases with the literature. Lastly, we consider some drawbacks of our implementation, one particular spurious feature of concern which developed with higher Weissenberg numbers, and describe some possible improvements.

---

[8] It is also appealing to remove $\hat{\mathbf{j}}$ entirely and replace it with the product of trace variables, $\hat{\mathbf{u}}\hat{\mathbf{T}}$. Here, $\hat{\mathbf{u}}$ would be the interface velocity and $\hat{\mathbf{T}}$ would be a new variable identified with the boundary restriction of $\mathbf{T}$. This choice was not made because it leads to an additional nonlinearity in the formulation.

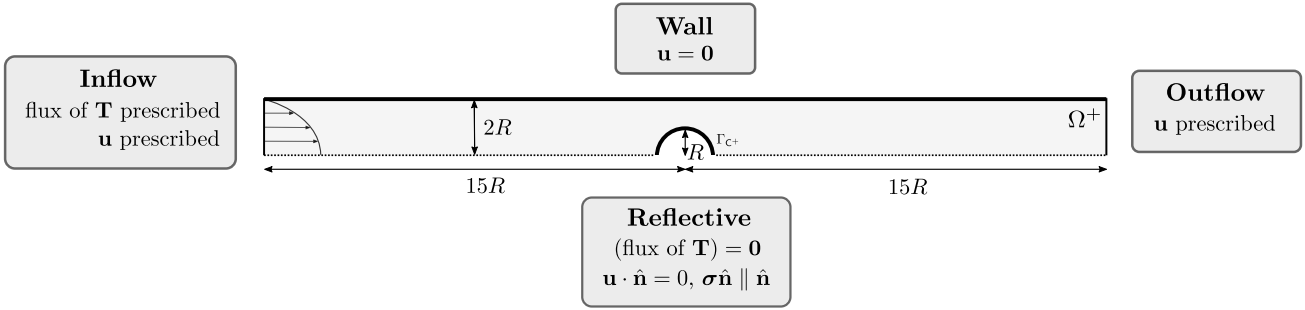

Figure 4.1: Diagram of flow domain with boundary names and the strong boundary conditions. Here, the prescribed flux of $\mathbf{T}$ is given in (4.2).

| Re | 0, 0.01, 0.1, 1 |
|---|---|
| Wi | 0.1, 0.2, 0.3, 0.4, 0.5, 0.6, 0.7, 0.8, 0.9, 1.0 |
| $\alpha$ | 0, 0.001, 0.01, 0.1 |
| $\beta$ | 0.59 |
| $p$ | 2 |
| $\mathrm{d}p$ | 2 |
| solver | MUMPS 5.0.1 [50, 51] |

Table 1: Parameters used in study.

*4.1. Set-up*

In this study, the Reynolds and Weissenberg numbers are defined as
$$\mathrm{Re} = \frac{\rho \bar{u} R}{\eta},$$
and
$$\mathrm{Wi} = \lambda \frac{\bar{u}}{R}.$$
Since we have used the radius of the cylinder as the length scale to define the Reynolds number for this problem, we choose to use it again in the definition of the graph norm, (3.6). That is, beginning at (3.5), we set $l_0 = R$.

As is standard in similar incompressible flow problems, we prescribed the inflow velocity to be that of the steady-state Poiseuille solution for flow in a channel, $\mathbf{u}(x,y) = \begin{pmatrix} u_1(y) \\ 0 \end{pmatrix}$, where $u_1(y) = \frac{3\bar{u}}{2}\left(1 - \frac{y^2}{4R^2}\right)$. Conveniently, this can be complemented with a simple solution for the extra-stress tensor, $\mathbf{T}(x,y) = \begin{pmatrix} T_{11}(y) & T_{12}(y) \\ T_{12}(y) & T_{22}(y) \end{pmatrix}$, where

$$T_{11}(y) = \frac{9\bar{u}^2 \lambda \eta_\mathrm{P} y^2}{8R^4}, \quad T_{12}(y) = \frac{-3\bar{u}\eta_\mathrm{P} y}{4R^2}, \quad T_{22}(y) = 0.$$

In some studies of this benchmark problem, the inflow value of the advected quantity, $\mathbf{T}$, is prescribed. Our formulation, however, does not readily allow this type of boundary condition.

We generally attempt to enforce boundary conditions, whenever possible, by explicitly constraining the interface variables *and not constraining the field variables*. However, the only interface variables which entered into the Lie derivative equations, (3.4), after integration by parts were the three independent components of the flux of the extra-stress, $\hat{\mathbf{j}}$. Prescribing the flux of $\mathbf{T}$ from the exact solution above, we set

$$\hat{\mathbf{j}} = \mathbf{T}\,\mathbf{u} \cdot \hat{\mathbf{n}} = -\mathbf{T}u_1 \tag{4.2}$$

at the inflow boundary.

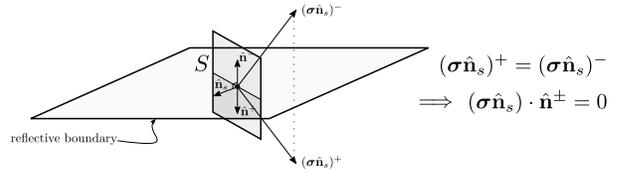

Figure 4.2: The illustration above is for a three-dimensional domain. Observe that the Cauchy stress vector must be continuous across the reflective boundary. However, it must also be symmetric under reflections passing the vector from $\Omega^+$ into $\Omega^-$. Therefore, at the reflective boundary, the stress vector must vanish in the direction normal to the interface boundary, $\boldsymbol{\sigma}\hat{\mathbf{n}}_s \perp \hat{\mathbf{n}}$. Since $\boldsymbol{\sigma}$ is symmetric and $\hat{\mathbf{n}}_s \perp \hat{\mathbf{n}}$ is arbitrary, all tangential components of $\hat{\mathbf{t}} = \boldsymbol{\sigma}\hat{\mathbf{n}}$ must vanish at the reflective boundary.

We chose the boundary conditions at the reflective boundary from physical intuition and motivation from the inflow boundary prescription. First of all, due to the symmetry of the solution, we anticipate zero *flux* of the extra stress tensor across the boundary. The logic for this is simple since any flux vector should exist in equal magnitude, but *reflected direction*, at the opposing point across the boundary. Following from this principle, the normal component of the flux must vanish at points where the flux is equal to its reflection, and so at all points on the reflective boundary, $\hat{\mathbf{j}} = \mathbf{0}$. Likewise, the velocity of the fluid normal to the boundary, the *mass flux*, must also vanish, $\mathbf{u} \cdot \hat{\mathbf{n}} = 0$. Indeed, we could conclude that the flux of $\mathbf{T}$ vanishes at the reflective boundary, plainly from this relationship on the fluid velocity, $\hat{\mathbf{j}} = \mathbf{T}\,\mathbf{u} \cdot \hat{\mathbf{n}} = \mathbf{0}$. The final boundary condition at the reflective boundary—the stress vector, $\hat{\mathbf{t}} = \boldsymbol{\sigma}\hat{\mathbf{n}}$, being parallel to the boundary normal—is not quite immediate although it can be argued similarly. Our argument for this condition is given in Figure 4.2.

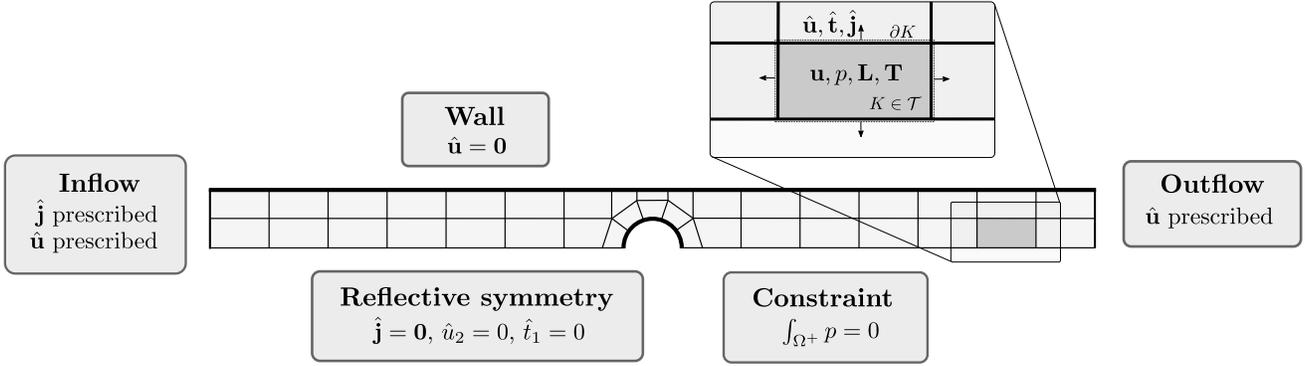

Figure 4.3: Initial mesh and DPG boundary conditions. Note that $\hat{\mathbf{u}} = \begin{pmatrix} \hat{u}_1 \\ \hat{u}_2 \end{pmatrix}$, $\hat{\mathbf{t}} = \begin{pmatrix} \hat{t}_1 \\ \hat{t}_2 \end{pmatrix}$, and $\hat{\mathbf{j}} = \begin{pmatrix} \hat{j}_{11} & \hat{j}_{12} \\ \hat{j}_{12} & \hat{j}_{22} \end{pmatrix}$.

At the outflow boundary, many different choices of boundary conditions are possible for this problem. We have chosen to present results from a fully prescribed outflow velocity field. Other choices such as zero outflow traction or mixed boundary conditions (zero tangential velocity and zero normal stress) were also tried by the authors; however, with the energy norm error indicator, both of these choices induced localized refinements near the outflow boundary. None of these other outflow boundary conditions had any discernible influence upon the computed drag coefficient values in our experiments. The presented boundary conditions—wherein the outflow velocity is prescribed from the Poiseuille solution—resulted in the smallest small-scale change of solution features near the outflow boundary and so also the fewest isolated refinements in that region.

Finally, boundary conditions at the walls of the channel and obstructing cylinder are simply taken to be of the standard no-slip type, $\mathbf{u} = \mathbf{0}$. Here, the flux of the extra stress tensor was not prescribed. Of course, prescribing the velocity field at all boundaries of the computational domain necessitated introducing a uniqueness constraint on the pressure of the system. In all experiments, we chose to enforce a zero-average pressure, $\int_{\Omega^+} p = 0$.

Given that the DPG method generally requires boundary conditions to be prescribed through the interface variables, we have reinterpreted them in Figure 4.3 as described above wherein we also depict the initial mesh used for each simulation. Notice that we always began our simulations with the same 36 elements and did not require any parameter continuation as did some other methods [12].

### 4.2. Creeping flow with the Oldroyd-B model

Here, we present the computed values of the drag coefficient for the Oldroyd-B model with Stokes flow coupling and the energy refinement strategy. In this situation, and throughout, all experiments performed with the energy strategy used the naive refinement marking scheme from [52]. That is, at each refinement step, for $\eta_{\max} = \max_{K \in \mathcal{T}} \eta_{V_K}$ and $\theta \in (0, 1)$, we marked elements for refinement if $\eta_{V_K} \geq \theta \eta_{\max}$. In our computations, this threshold parameter was consistently chosen to be $\theta = 0.2$.

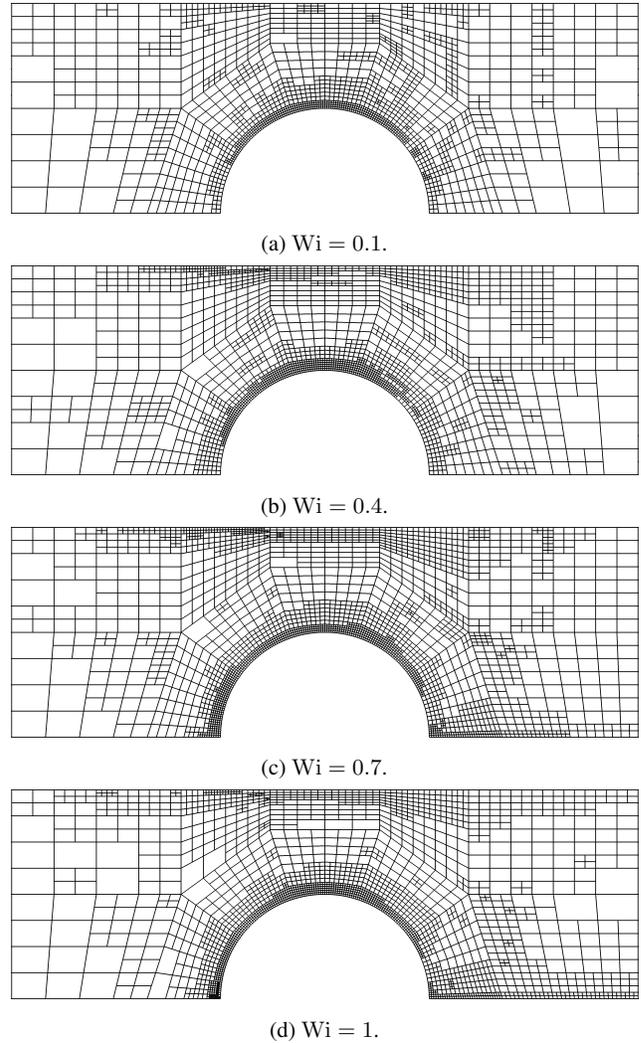

(a) Wi = 0.1.

(b) Wi = 0.4.

(c) Wi = 0.7.

(d) Wi = 1.

Figure 4.4: Close-up of meshes from the energy refinement strategy after five refinements for select values of the Weissenberg number.

|       | Wi = 0.1 |                  |         |                |            | Wi = 0.4 |                  |          |                |            |
| ----- | -------- | ---------------- | ------- | -------------- | ---------- | -------- | ---------------- | -------- | -------------- | ---------- |
| Ref. # | DoF     | Drag coefficient |         | Error estimate |            | DoF      | Drag coefficient |          | Error estimate |            |
|       |          | Field            | Flux    | Drag           | Energy     |          | Field            | Flux     | Drag           | Energy     |
| 0     | 5123     | 55.6760          | 56.1664 | 2.1454         | 0.3734     | 5123     | 33.9969          | 34.4629  | 1.7955         | 0.4296     |
| 1     | 19604    | 94.8348          | 95.0995 | 1.2605         | 0.2214     | 19604    | 91.1569          | 91.9478  | 2.0618         | 0.2358     |
| 2     | 24242    | 122.1375         | 122.1917| 0.5008         | 0.09841    | 30893    | 110.3167         | 110.6781 | 0.8045         | 0.09273    |
| 3     | 42488    | 129.4744         | 129.5010| 0.1501         | 0.03008    | 60839    | 116.6999         | 116.8479 | 0.3542         | 0.03740    |
| 4     | 94265    | 130.2408         | 130.2491| 0.04778        | 0.009435   | 108599   | 119.5058         | 119.5410 | 0.1114         | 0.01670    |
| 5     | 229358   | 130.3381         | 130.3420| 0.01543        | 0.003150   | 296321   | 120.4220         | 120.4297 | 0.03248        | 0.005749   |
| 6     | 571787   | 130.3566         | 130.3584| 0.006452       | 0.001031   | 648935   | 120.5616         | 120.5642 | 0.01102        | 0.002306   |
| 7     | 1122647  | 130.3608         | 130.3616| 0.003087       | 0.0004299  | 1265147  | 120.5818         | 120.5837 | 0.007196       | 0.001052   |
| 8     | 2477081  | 130.3620         | 130.3624| 0.001540       | 0.0001561  | 1923035  | 120.5879         | 120.5884 | 0.003450       | 0.0006169  |
| 9     | 4067018  | 130.3624         | 130.3626| 0.001950       | 0.00008547 | 3000239  | 120.5897         | 120.5902 | 0.002383       | 0.0003567  |
| 10    |          |                  |         |                |            | 3920696  | 120.5904         | 120.5906 | 0.001718       | 0.0002493  |

|       | Wi = 0.7 |                  |          |                |           | Wi = 1.0 |                  |          |                |           |
| ----- | -------- | ---------------- | -------- | -------------- | --------- | -------- | ---------------- | -------- | -------------- | --------- |
| Ref. # | DoF     | Drag coefficient |          | Error estimate |           | DoF      | Drag coefficient |          | Error estimate |           |
|       |          | Field            | Flux     | Drag           | Energy    |          | Field            | Flux     | Drag           | Energy    |
| 0     | 5123     | 31.5416          | 32.1447  | 1.9049         | 0.4374    | 5123     | 33.9760          | 34.6702  | 2.0224         | 0.4322    |
| 1     | 19604    | 85.4511          | 85.9790  | 1.3912         | 0.2057    | 19604    | 79.9033          | 80.3355  | 1.1545         | 0.2037    |
| 2     | 34076    | 105.0442         | 105.1755 | 0.4563         | 0.08371   | 36110    | 103.2292         | 103.2436 | 0.4292         | 0.09038   |
| 3     | 78146    | 112.5455         | 112.6402 | 0.2461         | 0.03462   | 80138    | 115.7134         | 115.7109 | 0.2457         | 0.04052   |
| 4     | 187970   | 115.0789         | 115.1318 | 0.1421         | 0.01771   | 204866   | 119.0960         | 119.1157 | 0.1127         | 0.02097   |
| 5     | 373802   | 116.6757         | 116.6810 | 0.04498        | 0.008825  | 353648   | 118.8900         | 118.8965 | 0.06567        | 0.02599   |
| 6     | 670955   | 116.9821         | 116.9929 | 0.03149        | 0.005488  | 362882   | 118.8117         | 118.8179 | 0.06386        | 0.01382   |
| 7     | 1137575  | 117.2007         | 117.2027 | 0.01587        | 0.003770  | 759218   | 118.3835         | 118.3891 | 0.04660        | 0.05420   |
| 8     | 1273412  | 117.1716         | 117.1739 | 0.01537        | 0.003152  | 834446   | 118.5343         | 118.5362 | 0.03567        | 0.009143  |
| 9     | 1683764  | 117.2265         | 117.2285 | 0.01433        | 0.002573  | 1169099  | 118.0705         | 118.0741 | 0.02763        | 0.007946  |
| 10    | 1785419  | 117.2284         | 117.2303 | 0.01437        | 0.002455  | 1427945  | 117.9592         | 117.9617 | 0.02655        | 0.005679  |
| 11    | 1860731  | 117.2319         | 117.2340 | 0.01438        | 0.002380  | 1643489  | 117.9875         | 117.9877 | 0.02284        | 0.004940  |
| 12    | 1904042  | 117.2345         | 117.2365 | 0.01433        | 0.002307  | 1955600  | 118.0815         | 118.0818 | 0.02638        | 0.004513  |
| 13    | 1955027  | 117.2379         | 117.2396 | 0.01367        | 0.002224  | 2406731  | 118.0868         | 118.0870 | 0.02369        | 0.005986  |
| 14    | 2007593  | 117.2395         | 117.2413 | 0.01369        | 0.002148  |          |                  |          |                |           |
| 15    | 2069720  | 117.2412         | 117.2430 | 0.01371        | 0.002066  |          |                  |          |                |           |
| 16    | 2134904  | 117.2435         | 117.2455 | 0.01517        | 0.001987  |          |                  |          |                |           |
| 17    | 2217206  | 117.2446         | 117.2466 | 0.01523        | 0.001775  |          |                  |          |                |           |
| 18    | 2929895  | 117.2890         | 117.2899 | 0.005264       | 0.001295  |          |                  |          |                |           |
| 19    | 3048377  | 117.2919         | 117.2928 | 0.005258       | 0.001240  |          |                  |          |                |           |
| 20    | 3130976  | 117.2925         | 117.2934 | 0.005254       | 0.001205  |          |                  |          |                |           |
| 21    | 3201977  | 117.2937         | 117.2946 | 0.005266       | 0.001166  |          |                  |          |                |           |
| 22    | 3262238  | 117.2942         | 117.2952 | 0.005278       | 0.001140  |          |                  |          |                |           |
| 23    | 3320939  | 117.2945         | 117.2955 | 0.005341       | 0.001115  |          |                  |          |                |           |
| 24    | 3379862  | 117.2951         | 117.2961 | 0.005332       | 0.001092  |          |                  |          |                |           |
| 25    | 3444287  | 117.2958         | 117.2968 | 0.005264       | 0.001070  |          |                  |          |                |           |
| 26    | 3539396  | 117.2971         | 117.2981 | 0.005215       | 0.001038  |          |                  |          |                |           |
| 27    | 3630017  | 117.2977         | 117.2987 | 0.005201       | 0.001010  |          |                  |          |                |           |
| 28    | 3717284  | 117.2984         | 117.2994 | 0.005165       | 0.0009821 |          |                  |          |                |           |
| 29    | 3806375  | 117.2990         | 117.3000 | 0.005107       | 0.0009549 |          |                  |          |                |           |
| 30    | 3893738  | 117.3001         | 117.3011 | 0.004956       | 0.0008700 |          |                  |          |                |           |

Table 2: Computed drag coefficient values are given above for Weissenberg numbers Wi = 0.1, 0.4, 0.7, and 1.0. We give the drag coefficient as computed from the flux variable $\hat{\mathbf{t}}$, as well as the drag coefficient as computed from the field variables. Reduction in the energy error as well as in the $L^2$ drag error is usually observed as each mesh is refined.

Table 2 demonstrates exactly how the numerical values of the drag coefficient change with the growing mesh for select values of the Weissenberg number. The computations ended when the MUMPS direct solver failed on our system. These failures could have been caused by insufficient computing resources for the problem size (an eventual barrier for all discretizations), ill conditioning, or simply ill-posedness and instability. Below, we attempt to determine, based on the data, the likely reason for each failure at each Wi. We also highlight the fact that the number of elements added to the mesh after each successive refinement could vary greatly with the Weissenberg number. With this in mind, even after just a handful of refinements, the mesh structure also varied greatly—as is demonstrated in Figure 4.4—with refinements in the wake of the cylinder becoming more pronounced with growing Wi.

We note the two different estimates of the drag coefficient presented in Table 2. The explanation for this is that by computing both the traction, $\hat{\mathbf{t}}$, as well as the field variable stress components, $p$, $\mathbf{L}$, and $\mathbf{T}$, we can make two *different* estimates of the drag coefficient:

$$\mathfrak{K}(\hat{\mathbf{t}}_h) = \frac{2}{\mu \bar{u}} \int_{\Gamma_{\text{c}+}} \hat{\mathbf{t}}_h \cdot \mathbf{e}_1 \,,$$

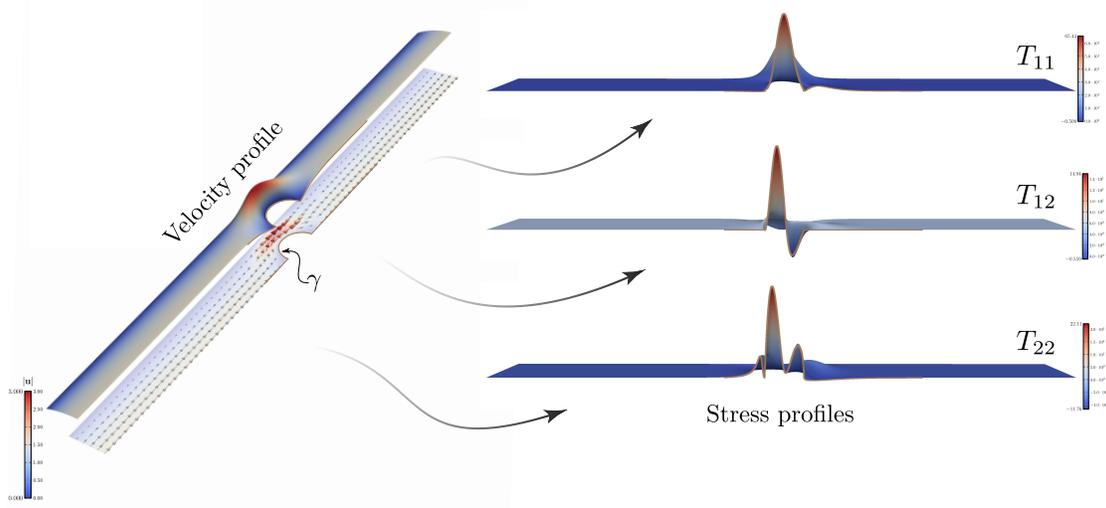

Figure 4.5: Extra-stress tensor components for Weissenberg number Wi $= 0.4$ in the Oldroyd-B fluid model. Computation from the $10^{\text{th}}$ refinement with the energy strategy (see Table 2). Displayed in the velocity profile is a surface plot of the velocity magnitude underlaid by a vector field plot showing the magnitude and direction of the velocity. The brown curve defines $\boldsymbol{\gamma}$, along which, the stress components are sampled in the coming figures.

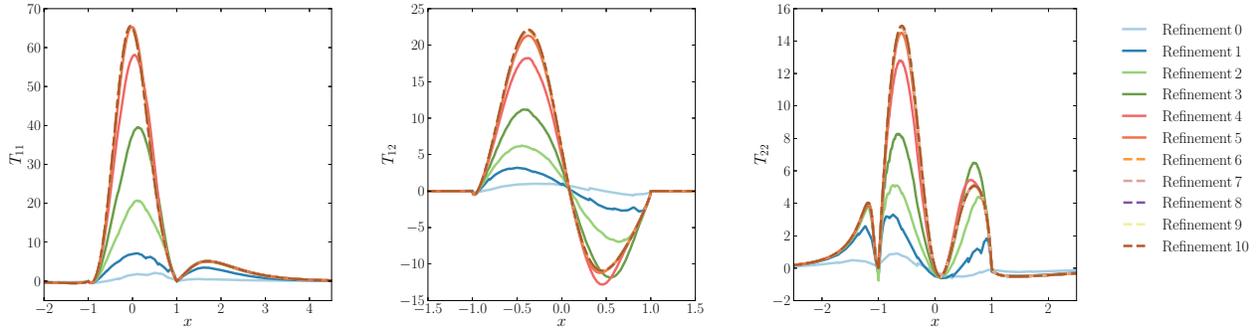

Figure 4.6: Convergence with mesh refinements in components of $\mathbf{T}$ along the curve $\boldsymbol{\gamma}$ (see Figure 4.5) for Wi $= 0.4$.

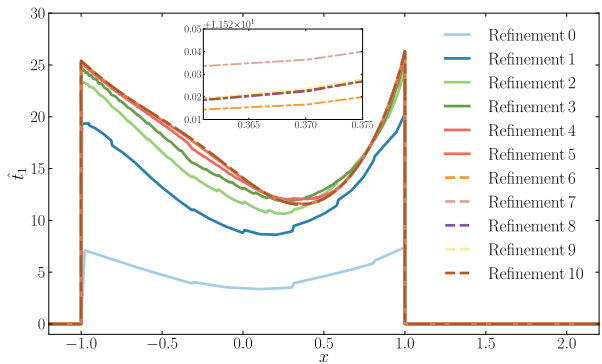

Figure 4.7: Convergence with mesh refinements in the first component of $\hat{\mathbf{t}}$ along $\boldsymbol{\gamma}$ for Wi $= 0.4$.

and
$$\mathfrak{K}(\boldsymbol{\sigma}_h \hat{\mathbf{n}}) = \frac{2}{\mu \bar{u}} \int_{\Gamma_{\text{C}+}} (\boldsymbol{\sigma}_h \hat{\mathbf{n}}) \cdot \mathbf{e}_1 \,,$$
where $\boldsymbol{\sigma}_h = -p_h \mathbf{I} + \eta_{\text{P}}(\mathbf{L}_h + \mathbf{L}_h^{\mathsf{T}}) + \mathbf{T}_h$. The first we call the *flux* estimate; and the second we call the *field* estimate.

Having two available estimates, we attempted to see if one was more accurate. Our inspection was both qualitative as well as quantitative. Initially, by plotting the profiles of the extra-stress components along the curve $\boldsymbol{\gamma}$ presented in Figure 4.5, we inspected the convergence with mesh refinement of $\mathbf{T}$ compared to $\hat{\mathbf{t}}$ for several values of the Weissenberg number. As illustrated in Figures 4.6 and 4.7, the profile of $t_1$ was generally less variable than any component of $\mathbf{T}$. This suggests that $\mathfrak{K}(\hat{\mathbf{t}}_h)$ could be more accurate simply because the accumulation of relative error in the field variables (to form $\boldsymbol{\sigma}_h$) is avoided. Another justification can be found by a simple examination of the theoretical energy spaces but we will not explore this here. It was, of course, the empirical evidence that was the strongest suggestion that $\mathfrak{K}(\hat{\mathbf{t}}_h)$ is the most accurate estimate. By observing the convergence behavior of both estimates and comparing them with the drag coefficient values reported in the literature, we found that the flux estimate was always slightly closer. In fact, although we are not sure of the reason, the field estimate was always slightly smaller than the flux estimate and, moreover, both approximations appeared to always converge to a steady value *from below*. We therefore only present the computed values of $\mathfrak{K}(\hat{\mathbf{t}}_h)$ in our results from now on.

Regardless of this determination, having two different estimates of the drag coefficient motivates a new *extrinsic* estimate of solution error. We define the drag error estimate,

$$\mathfrak{E}_\mathfrak{K} = |\Gamma_\mathsf{C}|^{1/2}\|(\hat{\mathbf{t}}_h - \boldsymbol{\sigma}_h\hat{\mathbf{n}})\cdot \mathbf{e}_1\|_{L^2(\Gamma_\mathsf{C})}.$$

We anticipate that for smooth enough solutions, this error will converge to zero. Figure 4.8 presents the behavior of this value as the mesh was refined for various values of Wi. Here, we see a relatively steady decrease in the drag error with refinement for $0.1 \leq$ Wi $\leq 0.4$, a lower rate and growing number of mesh refinements in the range $0.5 \leq$ Wi $\leq 0.7$, and progressively poorer results in the range $0.8 \leq$ Wi $\leq 1$.

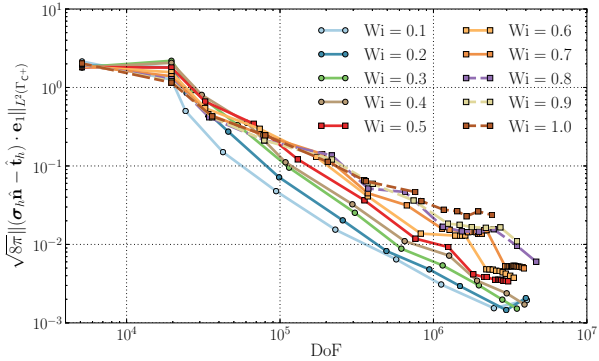

Figure 4.8: Convergence in drag error estimate versus degrees of freedom for the energy strategy. Markers indicate values at each refined mesh.

Understanding the behavior of the drag error with Weissenberg number was a concern in our study and for the range $0.5 \leq$ Wi $\leq 0.7$ some explanation can be given simply in the refinement pattern. Figure 4.9 demonstrates the growth of the mesh with the energy strategy for Wi $= 0.7$. Here, it is obvious that the energy error in the wake of the cylinder eventually dictates the mesh growth. For this reason, the mesh was rarely refined near the cylinder, so it is of little surprise that the drag error was not strongly affected after each refinement. It is well known that $T_{11}$ develops a strong internal layer in the wake of the confined cylinder as the Weissenberg number grows. This is depicted in Figure 4.10 and directly explains the refinement pattern we saw. Figure 4.11 demonstrates the evolution of this component of **T** along $\gamma$ as the mesh was refined when Wi $= 0.7$. For this Weissenberg number, this internal layer is known to be very difficult to reliably capture and our results, in accord with much of the literature, do not show convergence in the profile of $T_{11}$ in the wake of the cylinder.

As the drag coefficient is measured from solution values on the cylinder, we tested two ad-hoc refinement strategies which would refine the mesh more often near the cylinder and so hopefully increase the accuracy of our estimates. In strategy #1, we began with the same energy strategy as above, with the same threshold parameter, $\theta = 0.2$, except that at every step we also always refined each element with an edge lying on the cylinder boundary whether or not it was originally scheduled for refinement. In strategy #2, we similarly began with the energy error strategy except that at every step we enforced the refinement of each element with and edge lying within a distance of $0.1$ from the cylinder boundary. Close-ups of the sixth refined meshes for each of the three strategies when Wi $= 0.7$ are given in Figure 4.12.

Unfortunately, the ad-hoc strategies that we tested introduced issues of their own. In the first strategy, the relative scales of element sizes in the later meshes produced conditioning issues that led to failures in our solver. In the second strategy, the size of the narrow band about which we were enforcing mesh refinements was just large enough that all of our computations failed upon attempting the eighth refined mesh. In this second scenario, a slightly thinner band would likely have returned a more desirable final mesh; one with few enough degrees of freedom that our solver would not have crashed and we would have gotten a more accurate drag coefficient for our final data point. Determining the optimal band length was eventually abandoned as it was not in line with our research interests. Obviously, the energy strategy is not optimal for developing an accurate estimate of the drag coefficient. A goal-oriented approach would have been more desirable in this context, so we have begun developing such strategies with DPG for later work.

Table 3 compares the computed values of the drag coefficient with each of the three different adaptive strategies. Here, we see the best agreement with the literature in the energy strategy and so—while still recognizing its flaws—we decided to use it exclusively in our other studies.

A note must be made on our computations for Weissenberg numbers between $0.8$ and $1$. Although we will return to this again in Section 4.5, we mention that in this parameter interval, our nonlinear iterations failed to establish the expected quadratic rate of convergence in the Newton iterations as the mesh was being developed. Indeed, usually for coarse meshes, quadratic convergence of our Newton iterations was easily attained for all studied values of Wi. Only as the meshes were adaptively refined did the expected rate of converge falter for large Weissenberg numbers. This was true regardless of the refinement strategy we considered. Some reasons for this could possibly include degeneration of solution regularity but—considering similar results in the literature—most likely indicates that the problem we are solving is no longer well-posed. Indeed, some researchers have indicated that the solution becomes transient in this interval [47]. Others have indicated that the problem may be entirely ill-posed in this range due to the model allowing infinite extension of the viscoelastic fluid under finite elongation rates [9, 12]. Ultimately, we consider each result for the interval of Weissenberg numbers $0.8 \leq \text{Wi} \leq 1$—as well as for all our later results where quadratic convergence was not exhibited—dubious.

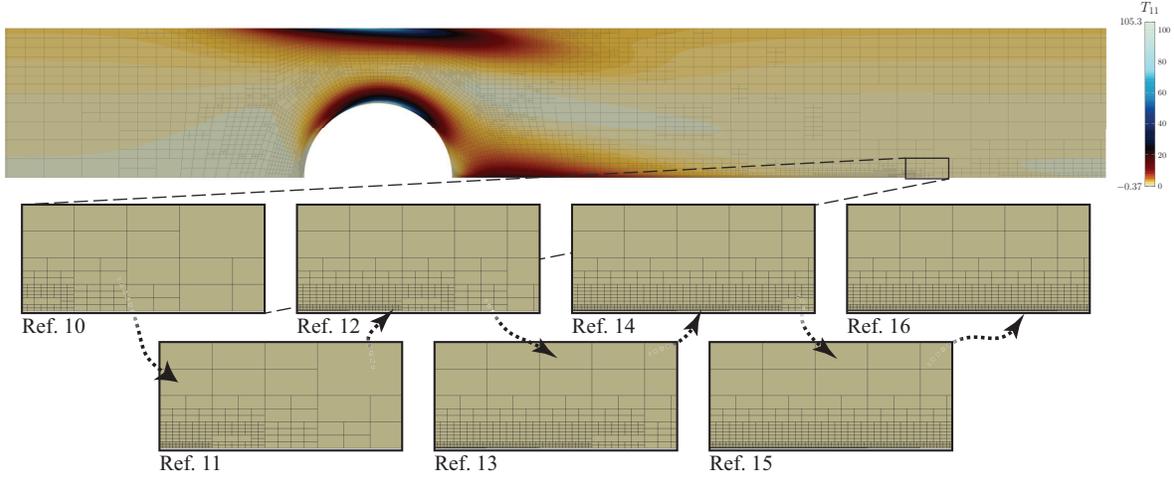

Figure 4.9: Close-up of a sequence of meshes from the Wi = 0.7 computation with the energy strategy. Notice the slow and sequential development of the mesh along the reflective boundary.

| Wi | Energy strategy | | Ad-hoc strategy #1 | | Ad-hoc strategy #2 | | Drag coefficient | | | | |
|---|---|---|---|---|---|---|---|---|---|---|---|
| | DoF(# Refs) | $\mathfrak{K}(\hat{\mathbf{t}}_h)$ | DoF(# Refs) | $\mathfrak{K}(\hat{\mathbf{t}}_h)$ | DoF(# Refs) | $\mathfrak{K}(\hat{\mathbf{t}}_h)$ | [12] | [47] | [9] | [48] | [49] |
| 0.1 | 4067018(9) | 130.3626 | 2630114(8) | 130.3625 | 2649878(7) | 130.3618 | 130.3626 | 130.364 | 130.363 | | 130.36 |
| 0.2 | 4001165(10) | 126.6251 | 3358379(9) | 126.6210 | 2504387(7) | 126.6241 | 126.6252 | 126.626 | 126.626 | | 126.62 |
| 0.3 | 3498518(10) | 123.1909 | 2214719(8) | 123.1904 | 2573072(7) | 123.1897 | 123.1912 | 123.192 | 123.193 | | 123.19 |
| 0.4 | 3920696(10) | 120.5906 | 2172425(8) | 120.5889 | 2777531(7) | 120.5885 | 120.5912 | 120.593 | 120.592 | | 120.59 |
| 0.5 | 3065843(17) | 118.8229 | 2093630(8) | 118.8150 | 2673770(7) | 118.8132 | 118.8260 | 118.826 | 118.836 | 118.827 | 118.83 |
| 0.6 | 3338165(19) | 117.7687 | 1858451(8) | 117.7370 | 2798306(7) | 117.7581 | 117.7752 | 117.776 | 117.775 | 117.775 | 117.77 |
| 0.7 | 3893738(30) | 117.3011 | 1649231(8*) | 117.1923 | 2606123(7*) | 117.2951 | 117.3157 | 117.316 | 117.315 | 117.291 | 117.32 |
| 0.8 | 4672934(12*) | 117.2973 | 1888487(8*) | 117.2091 | 2597321(7*) | 117.3057 | 117.3454 | 117.368 | 117.373 | 117.237 | 117.36 |
| 0.9 | 3503723(12*) | 117.5502 | 1847429(8*) | 117.5248 | 2607365(7*) | 117.6907 | 117.7678 | 117.812 | 117.787 | 117.503 | 117.79 |
| 1.0 | 2365391(14*) | 118.0873 | 930626(7*) | 118.7843 | 2506139(7*) | 118.5970 | | 118.550 | 118.501 | 118.030 | 118.49 |

Table 3: Comparison with results in literature with Stokes flow coupling. The superscript-* indicates that second order convergence was not reached in the final mesh.

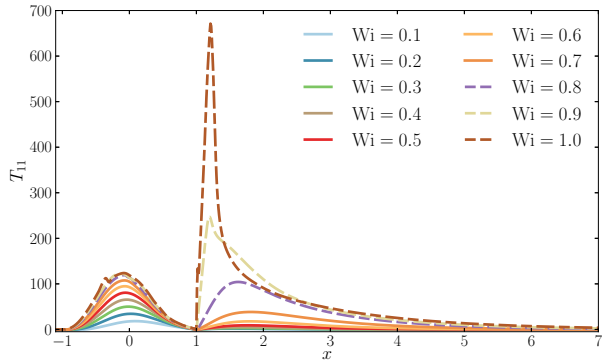

Figure 4.10: Profile curves for the $T_{11}$ component of the extra stress tensor along $\gamma$ among all Weissenberg numbers. Values taken from the energy strategy solutions at their final meshes (see Table 3).

Our observations in the loss of quadratic converge indicate an issue, not in the resolution of the problem and certainly not in the stability of the problem, but either in our discretization or in the underlying well-posedness itself. In order to gain further algorithmic insight, we have continued our analysis for the inertial Oldroyd-B model and the Giesekus model.

### 4.3. Effects of inertia in the Oldroyd-B model

In this subsection we investigate the effects of the advective term in the kinematic equations upon the Oldroyd-B model. Here, the drag coefficient was computed using the same energy strategy described above for Reynolds numbers Re = 0.01, 0.1, and 1. The results are collected together in Table 4. For a fixed Weissenberg number, the drag coefficient *grows* with the Reynolds number due to increasing velocity gradients.

Perhaps surprisingly, we see that our results match most closely with those in the literature as the Reynolds number is increased. The reason for this is suggested in Figure 4.13 where the profile of $T_{11}$ is plotted for each Reynolds number. Notably, as the Reynolds number grows, the values of $T_{11}$ decrease in the wake of the cylinder much faster than on the cylinder surface. This suggests that a smaller proportion of the elements will be marked for refinement downstream as Re grows, resulting in more accurate drag coefficient estimates. An inspection of the various meshes (not shown) verified this hypothesis.

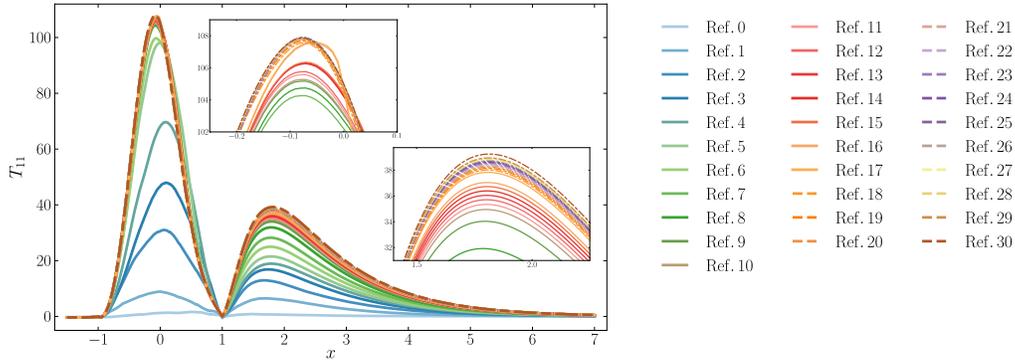

Figure 4.11: Convergence with mesh refinements in $T_{11}$ along $\gamma$ for Wi = 0.7.

| | Re = 0.01 | | | Re = 0.1 | | | Re = 1 | | |
|---|---|---|---|---|---|---|---|---|---|
| Wi | DoF(# Refs) | Drag coefficient | | Error est. | DoF(# Refs) | Drag coefficient | | Error est. | DoF(# Refs) | Drag coefficient | | Error est. |
| | | $\mathfrak{K}(\hat{\mathbf{t}}_h)$ | [47] | $\mathfrak{E}_\mathfrak{K}$ | | $\mathfrak{K}(\hat{\mathbf{t}}_h)$ | [47] | $\mathfrak{E}_\mathfrak{K}$ | | $\mathfrak{K}(\hat{\mathbf{t}}_h)$ | [47] | $\mathfrak{E}_\mathfrak{K}$ |
| 0.1 | 4065047(9) | 130.3628 | 130.364 | 0.001951 | 4208813(9) | 130.3667 | 130.368 | 0.001994 | 4146473(9) | 130.6075 | 130.609 | 0.002033 |
| 0.2 | 3993227(10) | 126.6259 | 126.627 | 0.002075 | 4635707(11) | 126.6352 | 126.636 | 0.002347 | 2849903(9) | 126.9366 | 126.938 | 0.001507 |
| 0.3 | 3489548(10) | 123.1926 | 123.194 | 0.001516 | 3691595(10) | 123.2095 | 123.211 | 0.001487 | 3760697(10) | 123.5957 | 123.597 | 0.001496 |
| 0.4 | 3933806(10) | 120.5933 | 120.595 | 0.001965 | 4666562(10) | 120.6197 | 120.622 | 0.001681 | 4485530(10) | 121.1038 | 121.106 | 0.001752 |
| 0.5 | 3097631(17) | 118.8269 | 118.831 | 0.003195 | 4472378(14) | 118.8654 | 118.868 | 0.002429 | 3573788(12) | 119.4567 | 119.460 | 0.002734 |
| 0.6 | 3718352(13) | 117.7752 | 117.781 | 0.003085 | 3535721(13) | 117.8234 | 117.831 | 0.003263 | 4364192(12) | 118.5380 | 118.542 | 0.002415 |
| 0.7 | 3900887(30) | 117.3079 | 117.323 | 0.004981 | 3924098(30) | 117.3717 | 117.387 | 0.005022 | 3982037(13) | 118.2221 | 118.233 | 0.004205 |
| 0.8 | 3503915(11*) | 117.2765 | 117.379 | 0.009565 | 3587420(14*) | 117.3533 | 117.459 | 0.01132 | 3306764(11) | 118.3971 | 118.455 | 0.009235 |
| 0.9 | 3560189(12*) | 117.5592 | 117.827 | 0.01100 | 3165626(22*) | 117.6836 | 117.925 | 0.01262 | 4106792(14) | 118.8741 | 119.096 | 0.01342 |
| 1.0 | 2708810(15*) | 118.1303 | 118.563 | 0.02817 | 2436611(13*) | 118.2237 | 118.697 | 0.02453 | 3063377(13) | 120.1455 | 120.057 | 0.02608 |

Table 4: Comparison with results in literature for the Oldroyd-B model with Navier-Stokes coupling. The superscript-* indicates that second order convergence was not reached in the final mesh. Notice that quadratic convergence was obtained for all Weissenberg numbers only in the case of the **largest** Reynolds number considered.

Figure 4.14 illustrates how the drag coefficient changes as the Reynolds number is increased. We see that the behavior is nearly identical until Re = 1.

### 4.4. Giesekus model

In this subsection we investigate the effects of the mobility factor in the Giesekus model with Stokes flow coupling. Here, the drag coefficient was computed using the same energy strategy described in Section 4.2 for $\alpha = 0.01, 0.1$, and 1. The results are collected together in Table 5. For a fixed Weissenberg number, the drag coefficient *decreases* as the mobility factor is increased. The decrease in the drag coefficient is due to the shear-thinning properties of the Giesekus fluid model.

Figure 4.15 depicts the profile of $T_{11}$ for each value of $\alpha$ and fixed Wi = 0.7. Notice that the scale of this variable is strongly dependent upon the order of magnitude of the mobility factor.

Figure 4.16 demonstrates how the drag coefficient pattern changed with Wi as the mobility factor was increased. Notably, the relationship became monotonic over the parameter range considered once $\alpha \geq 0.01$.

### 4.5. Further symmetry through penalty constraints

Reflect upon the results thus far. Generally, a reliable solution became more difficult to achieve as the Weissenberg number grew. This agrees with all literature on this problem and is clearly confirmed in Tables 3 through 5. Nevertheless, we wish to determine whether the computational challenges we have faced are likely to be caused by a particular feature of our method. The one possible issue which stands out to us is the choice of boundary conditions we have made.

Indeed, as remarked in Section 3.1, we created the interface variable $\hat{\mathbf{j}}$ so that the form $\hat{b}(\hat{u}, \mathfrak{v}^{\text{DPG}})$ would be bilinear. This caused us to enforce flux boundary conditions on the extra-stress tensor **T** on the boundary of the flow domain. A draw-back of this is that reflection symmetry in the velocity gradient and extra-stress tensor is not directly enforced. The purpose of this subsection is to examine this potential issue in isolation and determine if directly enforcing this additional symmetry would improve computation performance. We only consider the non-inertial Oldroyd-B model in this supplementary analysis.

Observe that $\frac{\partial u_2}{\partial x} = 0$ along the reflective boundary. Similarly, $\frac{\partial u_1}{\partial y} = 0$ from reflective symmetry of the velocity in the $x$-direction. Adding to these observations, the previous conclusion that $\sigma_{12} = 0$ along the reflective boundary, we conclude that, for any valid solution, $T_{12}$ must also vanish there. Analyzing the solutions of the Oldroyd-B model with varying Weissenberg numbers in Figure 4.17, we see that our computations fail to preserve this symmetry for each Wi $\geq$ 0.8. This suspicious circumstance calls each of these computations into question. The failure of our Newton iterations to converge quadratically for these parameter values may be related to degeneration of this symmetry condition, possibly indicating non-uniqueness of solutions.

| Wi | $\alpha = 0.001$ | | | $\alpha = 0.01$ | | | $\alpha = 0.1$ | | |
|---|---|---|---|---|---|---|---|---|---|
| | DoF(# Refs) | Drag coefficient | Error est. | DoF(# Refs) | Drag coefficient | Error est. | DoF(# Refs) | Drag coefficient | Error est. |
| | | $\mathfrak{K}(\hat{\mathbf{t}}_h)$ / [47] | $\mathfrak{E}_\mathfrak{K}$ | | $\mathfrak{K}(\hat{\mathbf{t}}_h)$ / [47] | $\mathfrak{E}_\mathfrak{K}$ | | $\mathfrak{K}(\hat{\mathbf{t}}_h)$ / [47] | $\mathfrak{E}_\mathfrak{K}$ |
| 0.1 | 4016738(9) | 130.2894 / 130.291 | 0.001946 | 3872162(9) | 129.6696 / 129.671 | 0.001885 | 4307282(12) | 125.5871 / 125.587 | 0.005766 |
| 0.2 | 3991931(10) | 126.3946 / 126.396 | 0.002062 | 4606553(10) | 124.6686 / 124.670 | 0.002196 | 4319063(12) | 117.1127 / 117.113 | 0.003802 |
| 0.3 | 3966476(9) | 122.7765 / 122.778 | 0.001459 | 3338810(10) | 120.0840 / 120.085 | 0.001569 | 4638866(13) | 111.0985 / 111.098 | 0.003748 |
| 0.4 | 4023257(10) | 119.9797 / 119.981 | 0.001604 | 3999086(10) | 116.5157 / 116.513 | 0.001584 | 4291865(13) | 106.8551 / 106.855 | 0.001851 |
| 0.5 | 3134828(18) | 118.0022 / 118.005 | 0.003174 | 3612272(10) | 113.8652 / 113.861 | 0.001756 | 4196786(13) | 103.7331 / 103.733 | 0.001864 |
| 0.6 | 2720891(12) | 116.7135 / 116.719 | 0.003383 | 4535513(17) | 111.9025 / 111.906 | 0.001768 | 4171346(12) | 101.3416 / 101.341 | 0.002014 |
| 0.7 | 4415051(12) | 115.9751 / 115.982 | 0.003502 | 4452977(13) | 110.4167 / 110.422 | 0.002077 | 4155383(12) | 99.4481 / 99.448 | 0.002250 |
| 0.8 | 3446018(13) | 115.6382 / 115.679 | 0.006507 | 4188863(14) | 109.2506 / 109.258 | 0.002519 | 4298489(12) | 97.9093 / 97.909 | 0.002235 |
| 0.9 | 4529801(13) | 115.6060 / 115.664 | 0.006488 | 3595223(13) | 108.2981 / 108.307 | 0.003393 | 4312274(14) | 96.6317 / 96.631 | 0.002221 |
| 1.0 | 3499793(14) | 115.7115 / 115.868 | 0.01087 | 4147868(15) | 107.4928 / 107.508 | 0.003247 | 4555955(12) | 95.5525 / 95.552 | 0.002206 |

Table 5: Comparison with results in literature for the Giesekus model with Stokes flow coupling.

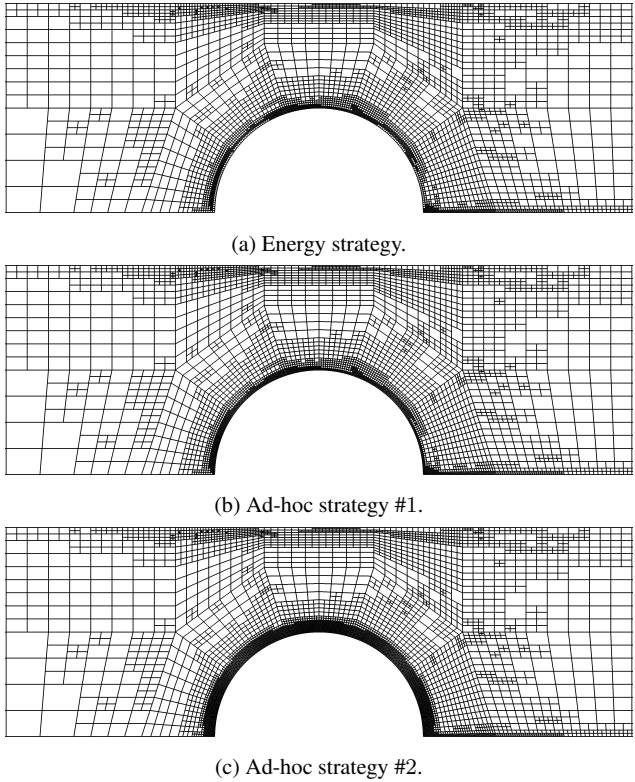

(a) Energy strategy.

(b) Ad-hoc strategy #1.

(c) Ad-hoc strategy #2.

Figure 4.12: Meshes near and in the wake of $\Gamma_{C+}$ from the three different strategies after six refinements. Wi = 0.7. At this point, the only large differences between the meshes are near the cylinder boundary.

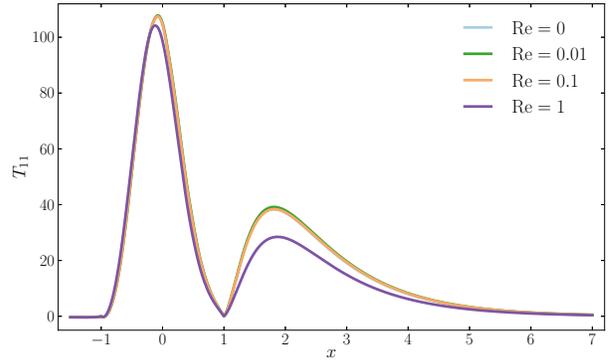

Figure 4.13: Profile of the $T_{11}$ component of the extra stress tensor for Wi = 0.7 along $\gamma$ with various Reynolds numbers.

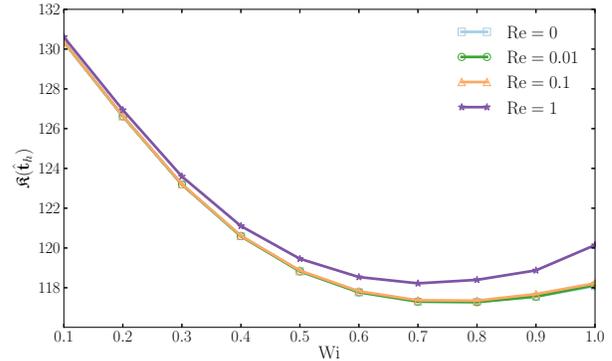

Figure 4.14: Dependence of the drag coefficient upon Reynolds number for Wi = 0.7.

As previously mentioned, given our variational formulation, we cannot influence the components of $T_{12}$ directly through boundary conditions. One possible remedy is to add penalty constraints on the field variables along the reflective boundary.

In our computations, symmetry in the $T_{12}$ variable was readily restored with the penalty constraint $T_{12}|_{\{y=0\}} = 0$ as depicted in Figure 4.18. Unfortunately, as is common with penalty methods, this constraint also increased the condition number of the linear systems. In Table 6 we have collected the results of our simulations and highlighted the scenarios in strongest disagreement with our previous computations. Notably, we see

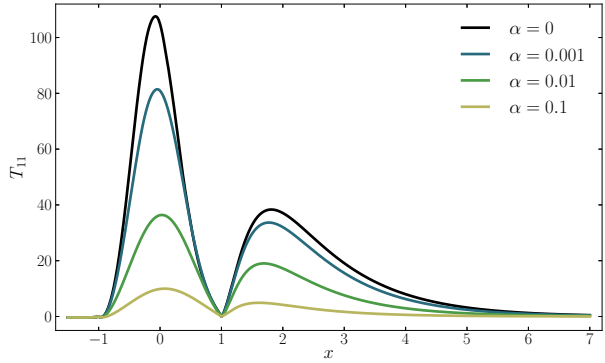

Figure 4.15: Profile of the $T_{11}$ component of the extra stress tensor for Wi = 0.7 along $\gamma$ with various values of $\alpha$.

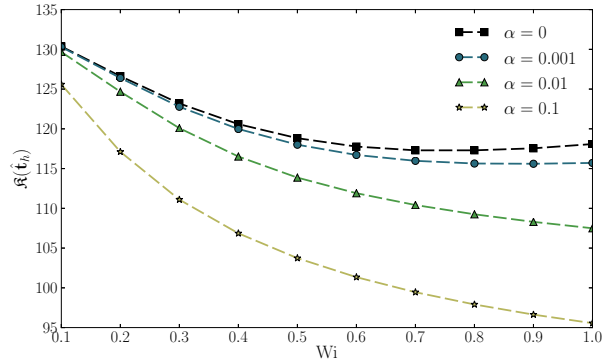

Figure 4.16: Dependence of the drag coefficient upon $\alpha$ for Wi = 0.7.

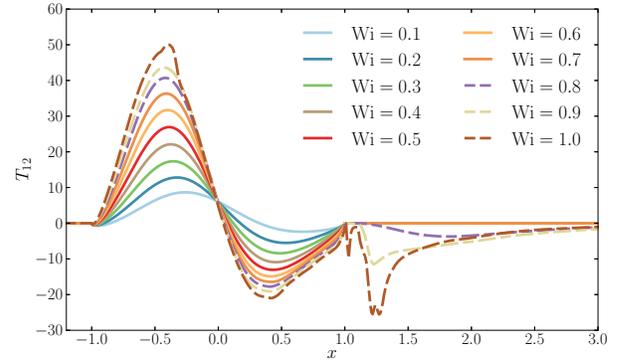

Figure 4.17: Profile curves for components of the extra stress tensor for various Weissenberg numbers from the final mesh created with the energy strategy. Observe the large inconsistency in the $T_{12}$ variable in the wake of the cylinder for the largest values of Wi. In these suspicious cases, this error was observed to compound through consecutive mesh refinements.

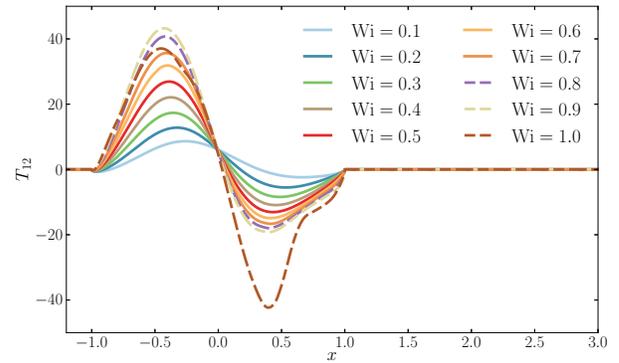

Figure 4.18: Profile curves for $T_{12}$ for various Weissenberg numbers from the final mesh created with the energy strategy **with penalty constraints**. Observe that the inconsistencies present in Figure 4.17 have been resolved in the wake. However, for most values of Wi, $T_{12}$ has not been affected on the cylinder. The curve for Wi = 1.0 is still almost certainly incorrect.

extremely good agreement for Wi ≤ 0.5 and the new estimate for Wi = 0.6 is slightly closer to those in the literature than before. Again, for Wi ≥ 0.8, we did not recover quadratic convergence in the Newton iterations, so we still consider these results dubious. This, therefore, does not justify further exploration into the inertial Oldroyd-B or non-inertial Giesekus models.

| Wi | # Refs | DoF | Drag coefficient, $\mathfrak{K}(\hat{\mathbf{t}}_h)$ | | Error est. $\mathfrak{E}_{\mathfrak{K}}$ |
|---|---|---|---|---|---|
| | | | penalty | previous | |
| 0.1 | 9/9 | 4069211 | 130.3626 | 130.3626 | 0.001950 |
| 0.2 | 10/10 | 4001165 | 126.6251 | 126.6251 | 0.002074 |
| 0.3 | 10/10 | 3492593 | 123.1909 | 123.1909 | 0.001516 |
| 0.4 | 10/10 | 3963839 | 120.5906 | 120.5906 | 0.002251 |
| 0.5 | 17/17 | 3072821 | 118.8230 | 118.8229 | 0.003432 |
| 0.6 | 14/19 | 4608965 | 117.7711 | 117.7687 | 0.003503 |
| 0.7 | 12/30 | 2496173 | 117.2786 | 117.3011 | 0.008573 |
| 0.8 | 12*/12* | 4742123 | 117.3268 | 117.2973 | 0.005994 |
| 0.9 | 13*/12* | 4579151 | 117.6461 | 117.5502 | 0.01077 |
| 1.0 | 6*/14* | 771287 | 120.1977 | 118.0873 | 0.06803 |

Table 6: Comparison with results in literature. Suspicious results are highlights in grey and superscript-* indicates that second order convergence was not reached in the final mesh.

## 5. Conclusion

In this paper we analyzed the steady flow of Oldroyd-B and Giesekus fluids around a confined cylinder with the DPG methodology. No stabilization was ever applied to our discretizations. We adaptively refined parameter-specific meshes for each simulation through a sequence of mesh refinements with hanging nodes marked by a local energy error estimate. This adaptive strategy was easily incorporated into our simulation and came naturally out of the DPG methodology. We used the Camellia software [3, 4] for all of our analysis.

Our qualitative comparisons with the existing literature were based on calculated estimates of the drag coefficient. Because our method contains two sets of solution variables (*field* and *interface*), we were able to estimate the drag coefficient in two different ways which led to an useful estimate of the solution error in the drag coefficient. We have called this error estimate, $\mathfrak{E}_\mathfrak{K}$, the (DPG) drag error estimate and we believe that it is unique in the literature. This error estimate is independent of the reported estimates in the literature and we anticipate that it may be useful for developing a stopping criterion in predictive modeling.

Unfortunately, the adaptive mesh refinement strategy we used in our studies was not constructed to produce optimal estimates of the drag coefficient. It is, instead, an intrinsic strategy that considers the solution error throughout the entire domain without bias. For this reason, for Weissenberg numbers generating strong boundary layers in the extra-stress tensor we observed diminishing agreement with the literature in terms of drag coefficient estimates. As the strength of the shock grew, we also saw increasing drag error estimates.

We briefly explored two ad-hoc refinement strategies to improve upon our drag coefficient estimates; however, our results in those examples portray the delicate nature of such strategies. Instead of advocating for a particular ad-hoc strategy for this problem, we have begun developing new goal-oriented refinement strategies for DPG methods which will appear in future research.

The ultimate conclusions of our simulations with respect to the well-posedness of this problem in certain Weissenberg number intervals are in general agreement with the literature. That is, we observed breakdown in our simulations for the Oldroyd-B model with Stokes flow coupling around $Wi = 0.7$. Our numerical experiments suggest to us that this is likely because of a lack of well-posedness in the problem itself for these values. In our simulations of the Oldroyd-B model with Navier-Stokes coupling, we saw the threshold Weissenberg number rise as the Reynolds number grew; similarly, the threshold Weissenberg number rose as the mobility factor was increased in the Giesekus model (with Stokes flow coupling). Our collections of drag coefficient estimates in these two simulations are in good agreement with [47] and improve as the Reynolds number or mobility factor is increased, respectively.

**Acknowledgements.** This work was partially supported with grants by NSF (DMS-1418822), AFOSR (FA9550-12-1-0484), ONR (N00014-15-1-2496) and by the German Research Foundation (DFG) under projects EL 741/3 and BE 3689/7.

## A. Appendix

Here we give a proof of the result used in Section 3.2 to explicitly derive the adjoint graph test norm that was used in our computations.

**Lemma A.1.** *Let $l \in X'$ where $X$ is Hilbert and $M \subset X$ is closed. Letting $N = M^\perp$, then*

$$\sup_{x_N \in N, x_M \in M} \frac{|l(x_N + x_M)|}{(\|x_N\|^2 + \|x_M\|^2)^{1/2}}$$
$$= \left( \sup_{x_N \in N} \frac{|l(x_N)|^2}{\|x_N\|^2} + \sup_{x_M \in M} \frac{|l(x_M)|^2}{\|x_M\|^2} \right)^{1/2}.$$

*Proof.* Recall the Riesz operator, $\mathcal{R}_X : X \to X'$, from (2.3). Let $x^\star = \mathcal{R}_X^{-1} l$ and so $(x^\star, \delta x)_X = l(\delta x)$ for all $\delta x \in X$ and so $\|x^\star\|_X = \|l\|_{X'}$. If we orthogonally decompose $x^\star = x_M^\star + x_N^\star$, $x_M^\star \perp x_N^\star$, then by orthogonality

$$(x_M^\star, \delta x_M)_X + (x_N^\star, \delta x_N)_X = l(\delta x_M) + l(\delta x_N),$$

for all $\delta x_M \in M, \delta x_N \in N$. Therefore, $x_M^\star = \mathcal{R}_M^{-1}(l|_M)$ and $x_N^\star = \mathcal{R}_N^{-1}(l|_N)$, and so

$$\|l\|_{X'}^2 = \|x^\star\|_X^2 = \|x_M^\star\|_M^2 + \|x_N^\star\|_N^2 = \|l|_M\|_{M'}^2 + \|l|_N\|_{N'}^2.$$

□